\documentclass[12pt]{article}
\usepackage{amsmath,amssymb,amsbsy,amsfonts,amsthm,latexsym,amsopn,amstext,amsxtra,euscript,amscd}
\usepackage{graphicx,shapepar}
\usepackage{amscd}	
\usepackage{latexsym}
\usepackage{epsfig}
\usepackage{arcs}

\title{Random triangles in planar regions containing a fixed point}

\author{Eugen J. Iona\c scu\\
\vspace{0.1cm} \\
\small Department of Mathematics, \\
\small Columbus State University\\
\small Columbus, GA 31907\\
\small Email: {\tt \{math\}@ejionascu.ro} }
\date{February $10^{th}$, 2018}
\def\n{\noindent}

\smallskip
\def\bsq{\blacksquare}
\def\eproof{$\hfill\bsq$\par}
 
\begin{document}
 \newtheorem{theorem}{\hspace{\parindent}
T{\scriptsize HEOREM}}[section]
\newtheorem{proposition}[theorem]
{\hspace{\parindent }P{\scriptsize ROPOSITION}}
\newtheorem{corollary}[theorem]
{\hspace{\parindent }C{\scriptsize OROLLARY}}
\newtheorem{lemma}[theorem]
{\hspace{\parindent }L{\scriptsize EMMA}}
\newtheorem{definition}[theorem]
{\hspace{\parindent }D{\scriptsize EFINITION}}
\newtheorem{problem}[theorem]
{\hspace{\parindent }P{\scriptsize ROBLEM}}
\newtheorem{conjecture}[theorem]
{\hspace{\parindent }C{\scriptsize ONJECTURE}}
\newtheorem{example}[theorem]
{\hspace{\parindent }E{\scriptsize XAMPLE}}
\newtheorem{remark}[theorem]
{\hspace{\parindent }R{\scriptsize EMARK}}
\renewcommand{\thetheorem}{\arabic{section}.\arabic{theorem}}
\renewcommand{\theenumi}{(\roman{enumi})}
\renewcommand{\labelenumi}{\theenumi}

\maketitle
\begin{abstract}
\smaller In this article we provide several exact formulae to calculate the probability that a random triangle chosen within a planar region (any Lebesgue measurable set of finite measure) contains a given fixed point $O$. These formulae are in terms of one integration of an appropriate function, with respect to a density function which depends of the point $O$.
The formulae provide another way to approach the Sylvester's Four-Point Problem as we show in the last section. A stability result is derived for the probability.
We recover the known probability in the case of an equilateral triangle and its center of mass: $\frac{2}{27}+20\frac{\ln 2}{81}$ (\cite{Kleitman}  and \cite{prekopa}).
We compute this probability in the case of a regular polygon and its center of mass for the point $O$.
Other families of regions are studied. For the family of Lima\c cons $r=a+\cos t$, $a>1$, and $O$ the origin of
the polar coordinates, the probability is  $\frac{1}{4}-\frac{12a^2(4a^2+1)}{(2a^2+1)^3\pi^2}$.
\end{abstract}

\section{Introduction}

Geometric probabilities are a rich source of beautiful examples in analysis when  calculations happen to add up in a neat way (see  \cite{Eisenberg&Sullivan}, \cite{GlenHall}, \cite{Ionascu&Prajitura}, and \cite{langford}).
This study is also such an instance. In what follows we are concerned with the question:
 {\it What is the probability that choosing three points at random $A$, $B$ and $C$ inside of a region $\bf R$ (uniform distribution with respect to the Euclidean planar area),
the triangle $\triangle ABC$ contains a fixed point $O$ which is in the interior of the region $\bf R$ or on its boundary?}

\[
\underset{Figure\ 1}{\underset{a}{\epsfig{file=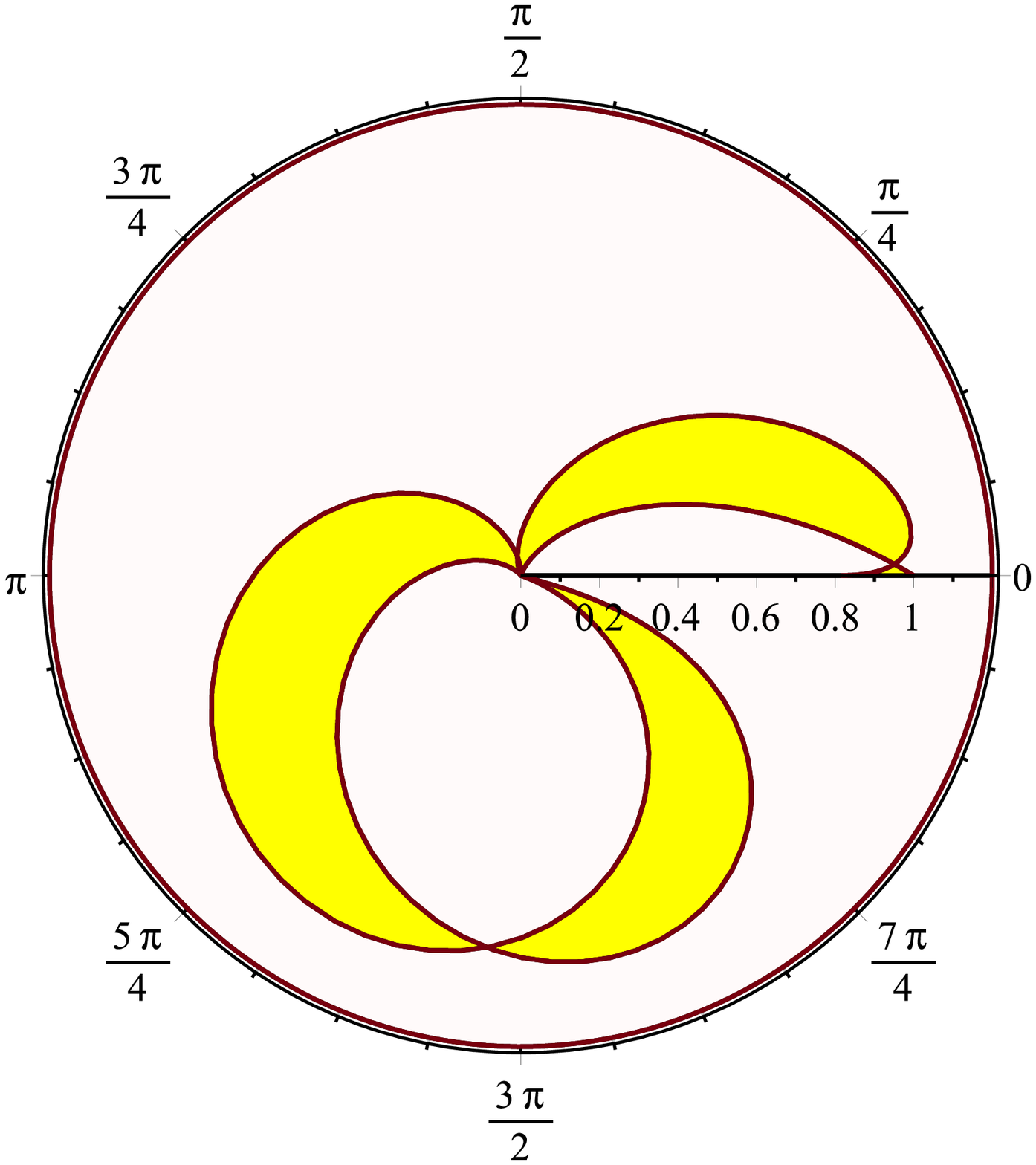,height=2.3in,width=2in}}\ \ \underset{b}{\epsfig{file=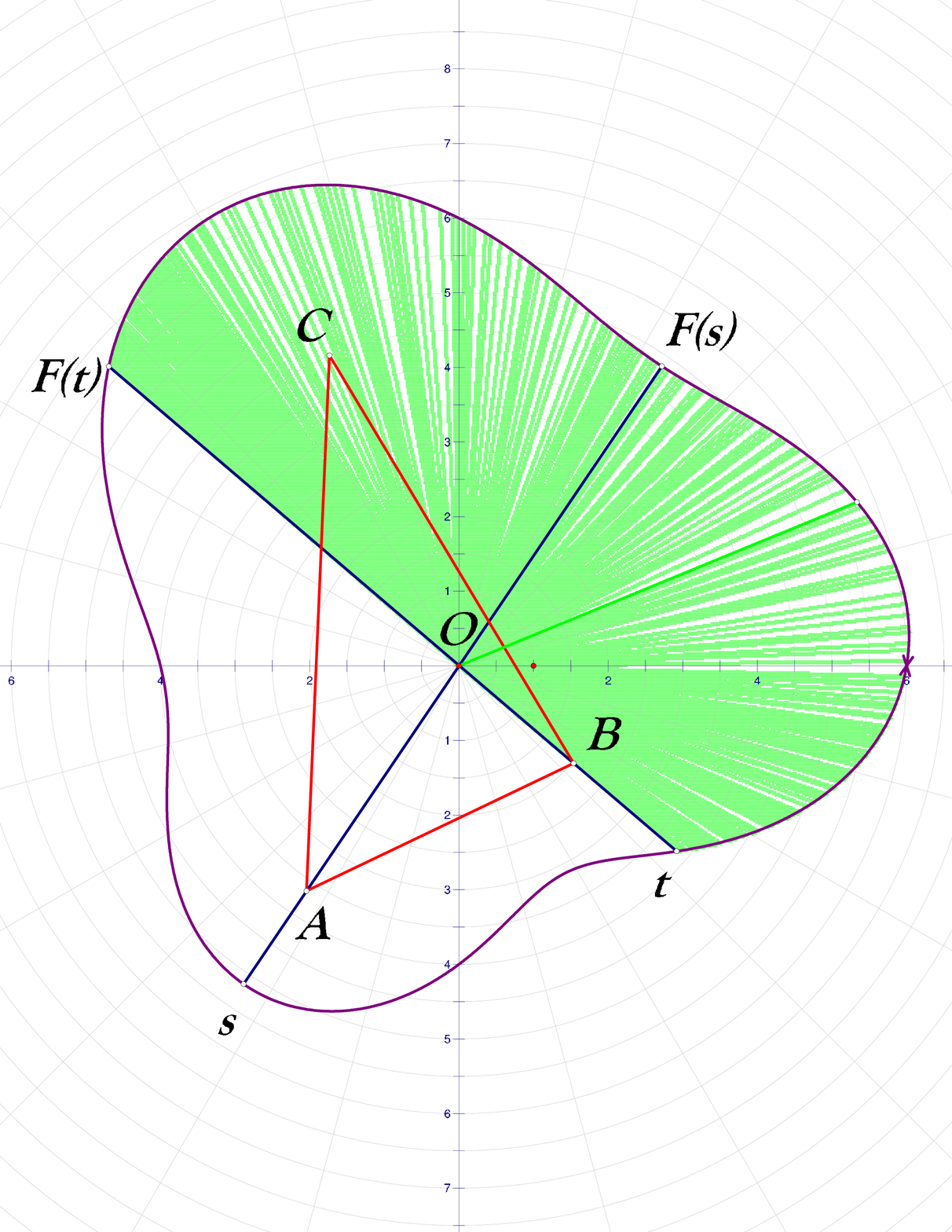,height=2.3in,width=2in}}}
\]

As a historical note we include what Adrian Baddeley had to say about this problem in \cite{AdrianB} (page 827):
``{\it An interesting problem concerns the probability $\rho(O)$ that a fixed point $O$ in an n-dimensional simplex $\bf S$ is covered by the convex hull of
$n-1$ independent and uniformly distributed random points of $\bf S$.  $\rm Andr\acute{a}s$ $\rm Pr\acute{e}kopa$ (\cite{prekopa}) conjectured that $\rho(O)$ is maximal when $O$ is the centroid (center of gravity) of $\bf S$; this has now been established for the planar triangle by $\it Hal\acute{a}sz$ and Kleitman (\cite{Kleitman}) using a refreshingly parsimonious approach.}" $\rm Andr\acute{a}s$ $\rm Pr\acute{e}kopa$ calculated this probability for an arbitrary triangle and $O\equiv G$ the center of mass and obtained $\frac{2}{27}+20\frac{\ln 2}{81}\approx 0.2452215$ (see \cite{prekopa}, page 21). Our approach is very close to this original derivation of this probability in the case of a triangular region.  A more recent paper of Norbert Henze (\cite{NorbertHenze}) studies a related problem (distribution of the area of a random triangle) in the context of the famous result of W. Blaschke which states that, among all regions of fixed area, the triangle and the circle provide bounds for the expected area of a random triangle.

\vspace{0.1in}

 For this introductory part we refer to Figure~1b. We make first a reduction to the boundary of the region, denoted here by $\partial \bf R$ which we will assume to be a smooth curve with the exception of a finite number of points. So, without loss of generality, we may assume that the area of the region $\bf R$ is equal to $1$, the point $O$ is the origin of a polar system of coordinates and  its boundary is given by a continuous, and piecewise differentiable function $r=r(\theta)$, $\theta\in [0,2\pi]$. Since the boundary is piecewise smooth we can use a parametrization of it with respect to its length, so any integral over the boundary included here, it is going to be with respect to its element of length.

If $A$ is an arbitrary point in $\bf R$, we let $s=s(A)$  its projection from $O$ on the boundary (Figure~1b).
 We consider a density distribution $f$ on the boundary, i.e., a continuous function $f:\partial {\bf R}\to [0,1]$ which satisfies:
\begin{equation}\label{distributionontheboundary}
\underset{\overarc{st}}{\int} f(x)dx=\frac{1}{2}\int_{\theta_s}^{\theta_t}r(\theta)^2d\theta,\ \theta_s\le \theta_t, \theta_t-\theta_s\le 2\pi,
\end{equation}

\n for every two points $s=(r(\theta_s),\theta_s)$ and $t=(r(\theta_t),\theta_t)$ (in polar)  on $\partial \bf R$, with the understanding that the arc $\overarc{st}$ is the oriented arc along the boundary
from $s$ to $t$ and the orientation is always counterclockwise.
By our assumption $\underset{\partial {\bf R}}{\int}f(x)dx=1$. We will see that our method works basically for any
parametrization of the boundary $s\to \gamma(s)$ and a density distribution $g$ such that for which $g(s)|\gamma'(s)|ds=\frac{1}{2}r(\theta_\gamma(s))^2d\theta_\gamma(s)$.

For a point $\bf t$ on the boundary let us denote by $F(\bf t)$ its projection through $O$ on the boundary: $F((r(\theta),\theta)):=
(r(\theta+\pi),\theta+\pi)$.  Clearly $F(F(t))=t$ for every $t\in \partial {\bf R}$, so $F$ is an involution on $\partial {\bf R}$.
Suppose we fix a point $\bf u$ on the boundary and introduce  two new functions: $$G(t)=\underset{\overarc{\bf ut}}{\int} f(x)dx\ \ \text{ and}\ H(t)=G(F(t))-G(t), \ \ t\in \overarc{\rm uv}$$ where $v=F(u)$.
The geometric interpretation of $H(t)$ is clearly the area of part of the region
 $\bf R$ which is to the right of the line $tF(t)$ (looking from $t$ to $F(t)$, see the green shaded area in Figure~1b).

The following formulae are of primary importance in this work.

\begin{theorem}\label{main} (i) The probability that choosing three points at random $A$, $B$ and $C$ inside of a region $\bf R$ of area $1$ (with uniform distribution with respect to the area),
the triangle $\triangle ABC$ contains the fixed interior point $O$ of the region $\bf R$ is given by
\begin{equation}\label{equation1}
 \boxed{{\cal P}=-G({\bf v})^2[3-2G({\bf v})]+6\underset{\ \ {\bf \overarc{\bf uv}}}{\int}H(t)\left[1-H(t)\right]f(t)dt,}
\end{equation} \n where ${\bf v}=F({\bf u})$.

(ii) Another expression for $\cal P$ is given by

\begin{equation}\label{equation11}
 \boxed{{\cal P}=\frac{1}{4}-\frac{(1-2x)^3}{4}-6\underset{\ \ {\bf \overarc{\bf uv}}}{\int}\left[\frac{1}{2}-H(t)\right]^2f(t)dt,}
\end{equation}
 \n where $x=G(v)$.

(iii) The point $\bf u\in \partial {\bf R}$ can be chosen in such a way the formula (\ref{equation1}) reduces to

\begin{equation}\label{equation2}
 \boxed{{\cal P}=\frac{1}{4}-6\underset{\ \ {\bf \overarc{\bf uv}}}{\int}\left[\frac{1}{2}-H(t)\right]^2f(t)dt.}
\end{equation}

 (iv) Independent of the point $u$, one can also write

\begin{equation}\label{equation3ind}
\boxed{{\cal P}=-\frac{1}{2}+3\underset{\ \ \bf \partial R}{\int}H(t)\left[1-H(t)\right]f(t)dt=\frac{1}{4}-3\underset{\bf \partial R}{\int}\left[\frac{1}{2}-H(t)\right]^2f(t)dt.}
\end{equation}

\end{theorem}

\n For polygonal regions the two functions $H$ and $f$ are not difficult to compute. For other types of regions one may want to use polar coordinates for the
computations of these functions or intrinsic parameterizations. For instance, if the curve is the lima\c con given by the polar equation $r=\frac{1}{3}\sqrt{\frac{2}{\pi}}(2+\cos \theta)$, $\theta\in [0,2\pi]$,
then the integral in (\ref{equation2}) reduces to  a pretty straightforward  calculation ${\cal P}=\frac{1}{4}-\frac{272}{243\pi^2} $.
Let us point out that one can arrive at the first formula in (\ref{equation3ind}) by a simple probabilistic argument as described in (\cite{Kleitman}).

\subsection{About the methods and the generality of the region $\bf R$}
We have written initially this paper under assumptions on $\bf R$ that are very restrictive. However, one can use the following idea in order to
relax the hypothesis on $\bf R$ to basically any Lebesgue measurable set in the plane of finite measure. Instead of working on the boundary of
$\bf R$, which in general can be as bad as the boundary of a fractal type set, we can work on a big circle of radius $r$ centered at $O$ (Figure~1a),
containing the region $\bf R$ if the region is bounded.
Even if $\bf R$ is not bounded,  the method can be used, but the formulae have to amended with a limit as $r\to \infty$.   The density function $f$ defined in (\ref{distributionontheboundary})
is going to be defined on the boundary of this circle, or equivalently on $[0,2\pi)$, by

\begin{equation}\label{distributionontheboundary2}
\underset{\overarc{st}}{\int} f(x)dx=\int_{\theta_s}^{\theta_t}f(\theta)d\theta=
 Area(\bf R \cap \angle {sOt}),\ \ \theta_s\le \theta_t, \theta_t-\theta_s\le \pi,
\end{equation}

Then all the formulae  (\ref{equation1})-(\ref{equation3ind}) can be recovered basically with the same proof.

\[
\underset{Figure\ 2, \ \ x\le x^2+y^2\le 1} { \epsfig{file=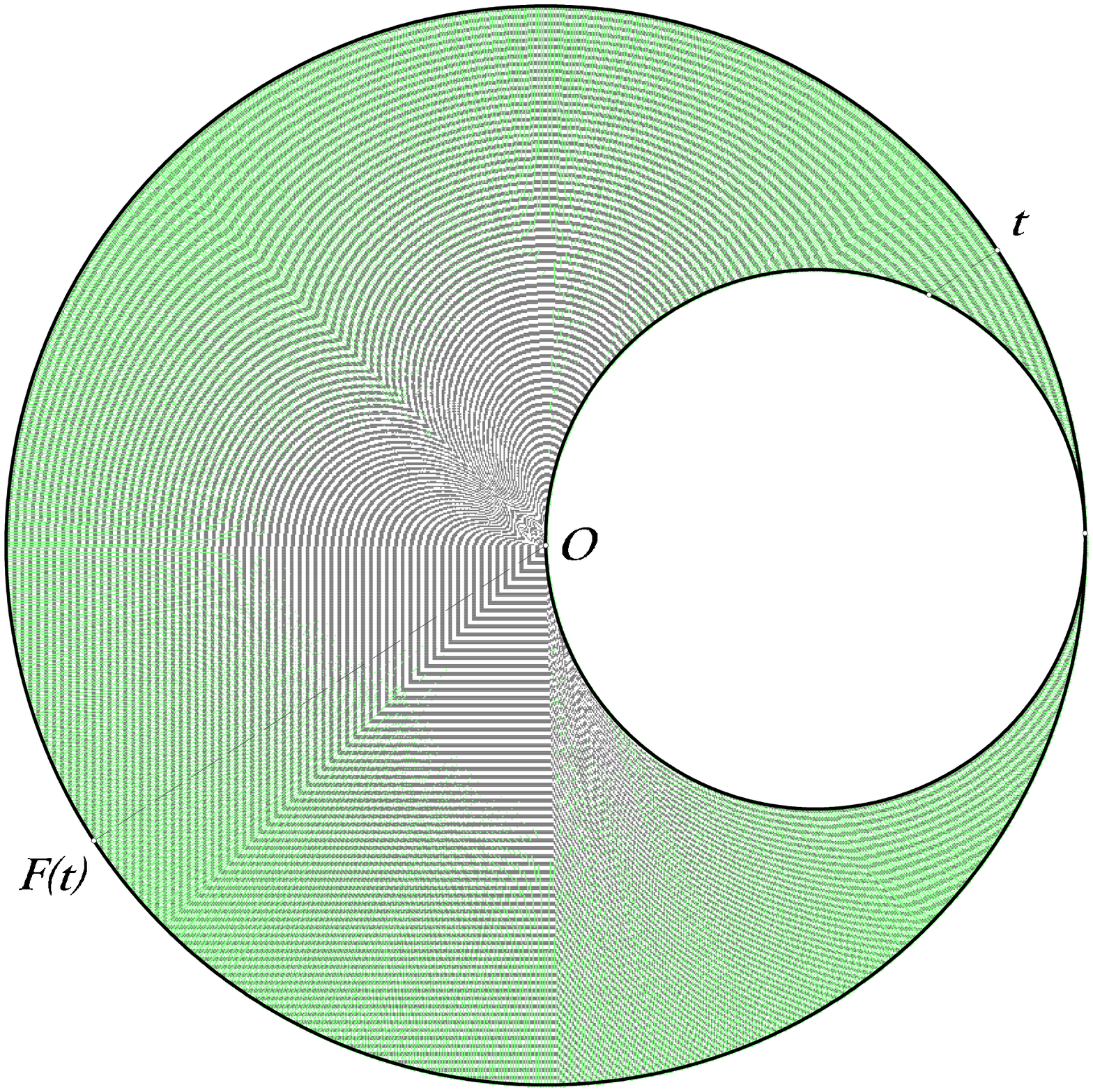,height=2.3in,width=2in} }
\]

For the region shown in Figure~2, with $O$ the origin of the rectangular coordinates,  we computed the functions $f$, $H$ and the value given by (\ref{equation2}), ${\cal P}=\frac{4\pi^2-5}{18\pi^2}\approx 0.1940774490$.
Computer simulations show that this is correct.

\section{Proof of Theorem~\ref{main}}

For two points on the boundary $s$ and $t$, we denote by ${\bf R}_{st}$ the region determined by intersection of $\bf R$ with the interior of the angle  $\angle F(s)OF(t)$. For every points $A$ and $B$ on the segments $\overline{Os}$ and $\overline{Ot}$, respectively,  then it is clear that
the third point $C$ must be in the region ${\cal R}_{st}$   in order for the point $O$ to be inside of $\triangle ABC$.
 This observation almost explains the next expression of the probability  we are interested in

 \begin{equation}\label{fimportanteq}
{\cal P}=\underset{\partial {\bf R}}{\int}\underset{\partial {\bf R}}{\int} Area({\cal R}_{st}) f(s)f(t)dsdt,
\end{equation}
\n but we will include a derivation of this formula. The probability we are looking for is equal to

$${\cal P}(O\in \triangle ABC)=\underset{ {\bf R}}{\int}\underset{ {\bf R}}{\int}\underset{ {\bf R}}{\int}\chi_{O\in \triangle ABC}(z_A,z_B,z_C)dz_Adz_Bdz_C.$$

\n where $dz_A$, $dz_B$ and $dz_C$ are the elements of area associated to the complex variable for each of the points $A$, $B$ and $C$ in ${\bf R}$ and
$\chi_E(\cdot)$ is the characteristic function of the set $E$, defined as usual by
$$\chi_E(\omega)=\begin{cases} 1\ \text{if}\ \omega \in E,\\ \\
0 \ \text{if}\ \omega \not \in E.
\end{cases}$$

\n Let us observe that, $\chi_{O\in \triangle ABC}(z_A,z_B,z_C)=\chi_{O\in \triangle s(A)BC}(z_A,z_B,z_C)$, which means that we can change to  polar coordinates in the integration over  ${\bf R}$ with respect to
the variable $z_A$:

 $$\underset{ {\bf R}^3}{\int} \chi_{O\in \triangle ABC}(z_A,z_B,z_C)dz_Adz_Bdz_C=
 \underset{ {\bf R}^2}{\int} \left[\int_0^{2\pi}\int_0^{r(\theta)}\chi_{O\in \triangle s(A)BC}(re^{i\theta},z_B,z_C)rdrd\theta \right]dz_Bdz_C.$$

\n In the integrant above $\chi_{O\in \triangle s(A)BC}(re^{i\theta},z_B,z_C)=\chi_{O\in \triangle s(A)BC}(s(A),z_B,z_C)$ which means this is a constant function along the segment
$Os(A)$. This results in the integration with respect to $r$ to be $\frac{r(\theta)^2}{2}$. Using the definition of the density distribution $f$ along the boundary, we can  change the  variable
from $\theta$ to $s$,   $\frac{r(\theta)^2}{2}d\theta=f(s)ds$, and obtain that

$$ {\cal P}= \underset{ {\bf R}^2}{\int} \left[\underset{\partial {\bf R}}{\int} \chi_{O\in \triangle s BC}(s,z_B,z_C)f(s)ds \right]dz_Bdz_C.$$

\n Doing the same thing for the variable $B$, we end up with

$$ {\cal P}= \underset{ {\bf R}}{\int}\left[\underset{\partial {\bf R}}{\int} \underset{\partial {\bf R}}{\int} \chi_{O\in \triangle s tC}(s,t,z_C)f(s)f(t)dtds \right]dz_C.$$

\n On the other hand, for $s$ and  $t$ fixed on the boundary, $$\underset{ {\bf R}}{\int}\chi_{O\in \triangle s tC}(s,t,z_C)dz_C=Area({\cal R}_{st})$$
\n so, interchanging the integrals and integrating with respect to $z_C$ first, we obtain (\ref{fimportanteq}).

\vspace{0.1in}

From now on, we are going to manipulate this to bring it in the form of (\ref{equation1}).  We could have started with this representation, which is basically the first intuitive reduction
and the purpose of introducing the distribution $f$ on the boundary.  For a $t\in \partial {\bf R}$ fixed, and $s\in tF(t)$,
the area of ${\cal R}_{st}$ is given by $\underset{F(t)F(s)}{\int}f(x)dx$. If $s\in F(t)t$,
the area of ${\cal R}_{st}$ is given by $\underset{tF(s)}{\int}f(x)dx$ (we will skip the arc notation for simplicity).

\n In (\ref{fimportanteq}), we split the boundary into two parts $\partial {\bf R}=\overarc{\bf uv} \cup \overarc{\bf vu}$:
and the corresponding integrals ${\cal P}=I_1+I_2$.

Then, the inner integration we split into four parts (Figure~2),

$$\begin{array}{c}
I_1=\underset{\bf uv}{\int}[\underset{\bf ut}{\int} (G(F(t))-G(F(s))f(s)ds+\underset{\bf tv}{\int}(G(F(s))-G(F(t)))f(s)ds+\\ \\
+\underset{\bf vF(t)}{\int}(1+G(F(s))-G(F(t)))f(s)ds+\underset{\bf F(t)F(v)}{\int}(G(F(t))-G(F(s)))f(s)ds]f(t)dt=J_1+J_2.
\end{array}$$
\n The last integration was split into two integrals as follows

$$\begin{array}{c}
J_1=\underset{\bf uv}{\int}[-\underset{\bf ut}{\int}G(F(s))f(s)ds+\underset{\bf tv}{\int}G(F(s))f(s)ds+\\ \\
+\underset{\bf vF(t)}{\int}G(F(s))f(s)ds-\underset{\bf F(t)F(v)}{\int}G(F(s)))f(s)ds]f(t)dt,\ \ \text{and}
\end{array}$$

$$\begin{array}{c}
J_2=\underset{\bf uv}{\int}[\underset{\bf ut}{\int} G(F(t))f(s)ds-\underset{\bf tv}{\int}G(F(t))f(s)ds+\\ \\
+\underset{\bf vF(t)}{\int}(1-G(F(t)))f(s)ds+\underset{\bf F(t)F(v)}{\int}G(F(t))f(s)ds]f(t)dt.
\end{array}$$
The second integral $J_2$ can be simplified to

$$\begin{array}{c}
J_2=\underset{\bf uv}{\int}[G(F(t))G(t)ds-G(F(t))(G({\bf v})-G(t))+\\ \\
+(1-G(F(t)))(G(F(t))-G({\bf v}))+G(F(t))(1-G(F(t))]f(t)dt=\\ \\
=\underset{\bf uv}{\int}[2G(F(t))(1-G(F(t))+G(t))-G({\bf v})]f(t)dt=\\ \\
=\underset{\bf uv}{\int}2G(F(t))\left[1-G(F(t))+G(t)\right]f(t)dt-G({\bf v})^2.
\end{array}$$

\[
\underset{Figure\ 3,\ \text{Graph of F }}{\epsfig{file=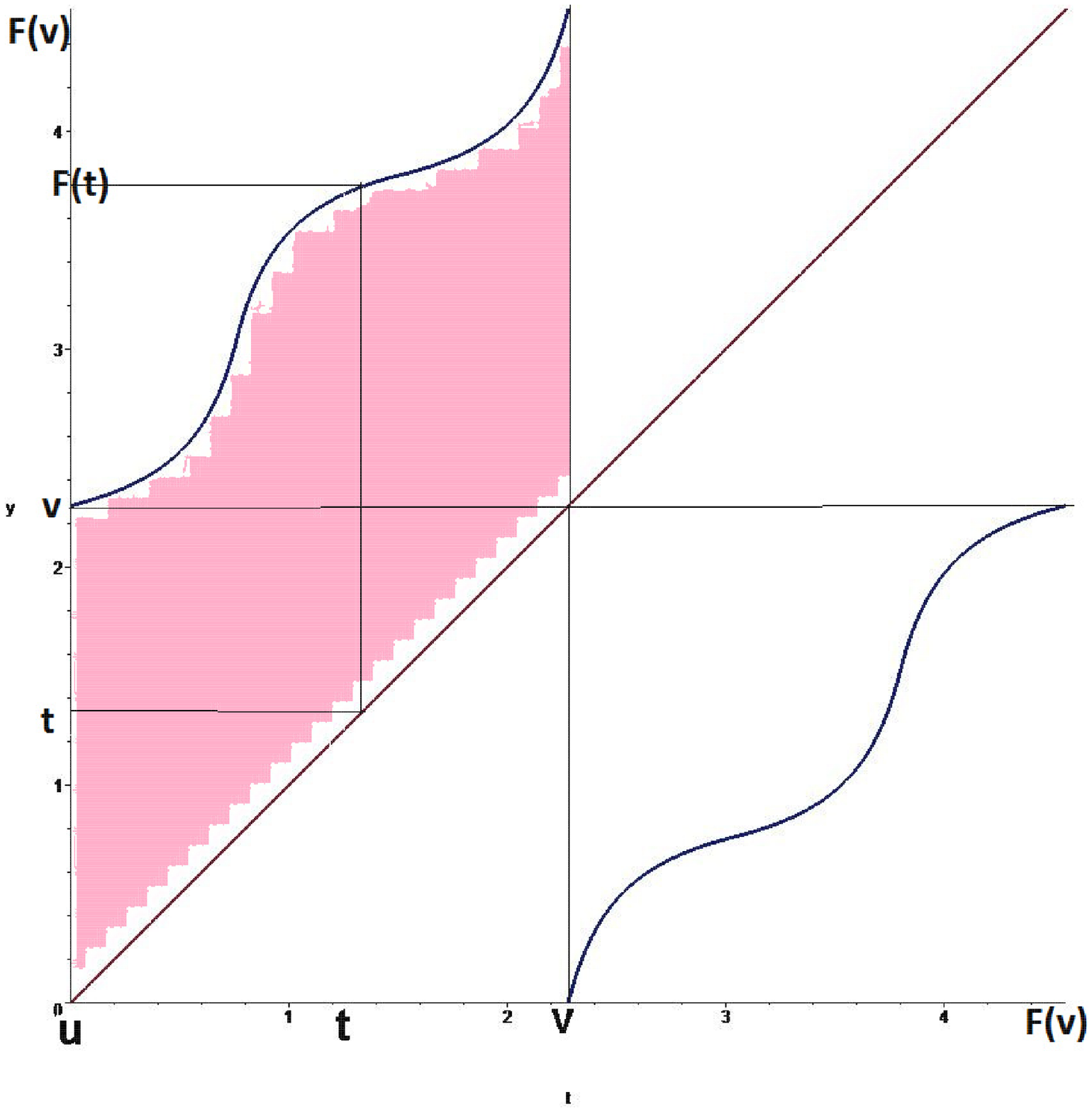,height=3.3in,width=3in,angle=-90}}
\]
\n For $J_1$ we will interchange integration order of integration in each one of those four integrals (Figure~3), for instance,

$$\begin{array}{c} \underset{\bf uv}{\int}\underset{\bf ut}{\int}G(F(s))f(s)f(t)dsdt=\underset{\bf uv}{\int}\underset{\bf sv}{\int}G(F(s))f(t)f(s)dtds=\\ \\
\underset{\bf uv}{\int} G(F(s))[G({\bf v})-G(s)]f(s)ds.
\end{array}$$
Similarly, we get
$$\begin{array}{c} \underset{\bf uv}{\int}\underset{\bf tv}{\int}G(F(s))f(s)f(t)dsdt=\underset{\bf uv}{\int}\underset{\bf us}{\int}G(F(s))f(t)f(s)dtds=\\ \\
\underset{\bf uv}{\int} G(F(s))G(s)f(s)ds,
\end{array}$$

$$\begin{array}{c} \underset{\bf uv}{\int}\underset{\bf vF(t)}{\int}G(F(s))f(s)f(t)dsdt=\underset{\bf vu}{\int}\underset{\bf F(s)v}{\int}G(F(s))f(t)f(s)dtds=\\ \\
\underset{\bf vu}{\int} G(F(s))[G({\bf v})-G(F(s)]f(s)ds,
\end{array}$$

\n and finally

$$\begin{array}{c} \underset{\bf uv}{\int}\underset{\bf F(t)F({\bf v})}{\int}G(F(s))f(s)f(t)dsdt=\underset{\bf vu}{\int}\underset{\bf uF(s)}{\int}G(F(s))f(t)f(s)dtds=\\ \\
\underset{\bf vu}{\int} G(F(s))G(F(s)f(s)ds.
\end{array}$$

\n Putting all these last four integrals together,  gives
$$J_1=\underset{\bf uv}{\int} G(F(s))[2G(s)-G({\bf v})]f(s)ds+\underset{\bf vu}{\int} G(F(s))[G({\bf v})-2G(F(s)]f(s)ds.$$

\n Hence, adding $J_1$ and $J_2$ one obtains

\begin{equation}\label{eq1}
\begin{array}{c} I_1=\underset{\bf uv}{\int}G(F(t))\left[2-2G(F(t))+4G(t)-G({\bf v})\right]f(t)dt+\\ \\ \underset{\bf vu}{\int} G(F(t))[G({\bf v})-2G(F(t)]f(t)dt-G({\bf v})^2.
\end{array}\end{equation}

\n Now, we handle $I_2$ the same way

$$\begin{array}{c}
I_2=\underset{\bf vu}{\int}[\underset{\bf uF(t)}{\int} (G(F(s))-G(F(t))f(s)ds+\underset{F(t){\bf v}}{\int}(1+G(F(t))-G(F(s)))f(s)ds+\\ \\
+\underset{\bf vt}{\int}(G(F(t))-G(F(s)))f(s)ds+\underset{\bf tF(v)}{\int}(G(F(s))-G(F(t)))f(s)ds]f(t)dt=K_1+K_2.
\end{array}$$

\n As before, we define the terms $K_1$ and $K_2$ in the following way:

$$\begin{array}{c}
K_1=\underset{\bf vu}{\int}[\underset{\bf uF(t)}{\int} G(F(s))f(s)ds-\underset{\bf F(t)v}{\int}G(F(s))f(s)ds+\\ \\
-\underset{\bf vt}{\int}G(F(s))f(s)ds+\underset{\bf tF(v)}{\int}G(F(s))f(s)ds]f(t)dt,
\end{array}$$

$$\begin{array}{c}
K_2=\underset{\bf vu}{\int}[\underset{UF(t)}{\int} -G(F(t)f(s)ds+\underset{\bf F(t)v}{\int}(1+G(F(t)))f(s)ds+\\ \\
+\underset{\bf vt}{\int}G(F(t))f(s)ds-\underset{\bf tF(v)}{\int}G(F(t))f(s)ds]f(t)dt.
\end{array}$$

\n Let's start with the easy part $K_2$, and observe that
$$\begin{array}{c}
K_2=\underset{\bf vu}{\int}[ -G(F(t))G(F(t))+(1+G(F(t))(G({\bf v})-G(F(t)))+\\ \\
+G(F(t)(G(t)-G({\bf v}))-G(F(t))(1-G(t))]f(t)dt \Rightarrow
\end{array}$$

$$\begin{array}{c}
K_2=\underset{\bf vu}{\int}[ -2G(F(t))[1-G(t)+G(F(t))]f(t)dt+G({\bf v})(1-G({\bf v})).
\end{array}$$

\n For the first term in $K_1$ we have

$$\begin{array}{c} \underset{\bf vu}{\int}\underset{\bf uF(t)}{\int} G(F(s))f(s)f(t)dsdt=\underset{\bf uv}{\int}\underset{\bf F(s)u}{\int} G(F(s))f(t)f(s)dtds=\\ \\
=\underset{\bf uv}{\int}G(F(s))(1-G(F(s))f(s)ds.
\end{array}$$

\n For the second term, we get

$$\begin{array}{c} \underset{\bf vu}{\int}\underset{\bf F(t)v}{\int} G(F(s))f(s)f(t)dsdt=\underset{\bf uv}{\int}\underset{\bf vF(s)}{\int} G(F(s))f(t)f(s)dtds=\\ \\
=\underset{\bf uv}{\int}G(F(s))(G(F(s)-G({\bf v}))f(s)ds.
\end{array}$$

\n Then, the third term

$$\begin{array}{c} \underset{\bf vu}{\int}\underset{\bf vt}{\int} G(F(s))f(s)f(t)dsdt=\underset{\bf vu}{\int}\underset{tU}{\int} G(F(s))f(t)f(s)dtds=\\ \\
=\underset{\bf vu}{\int}G(F(s))(1-G(t))f(s)ds,
\end{array}$$
\n and finally

$$\begin{array}{c} \underset{\bf vu}{\int}\underset{\bf tu}{\int} G(F(s))f(s)f(t)dsdt=\underset{\bf vu}{\int}\underset{\bf vt}{\int} G(F(s))f(t)f(s)dtds=\\ \\
=\underset{\bf vu}{\int}G(F(s))(G(t)-G({\bf v}))f(s)ds.
\end{array}$$

\n Hence, all these together gives

$$K_1=\underset{\bf uv}{\int}G(F(s))(1-2G(F(s))+G({\bf v}))f(s)ds+\underset{\bf vu}{\int}G(F(s))(2G(t)-1-G({\bf v}))f(s)ds.$$

This gives the second integral $I_2$:

$$\begin{array}{c} I_2=\underset{\bf uv}{\int}G(F(s))(1-2G(F(s))+G({\bf v}))f(s)ds+\\ \\
\underset{\bf vu}{\int}G(F(s))(4G(t)-3-2G(F(t))-G({\bf v}))f(s)ds+G({\bf v})(1-G({\bf v})),\end{array}$$

\n and so

\begin{equation}\label{eq3}
\boxed{\begin{array}{c}
{\cal P}=\underset{\bf uv}{\int}G(F(t))\left[3-4G(F(t))+4G(t)\right]f(t)dt-\\ \\
-\underset{\bf vu}{\int}G(F(t))\left[3+4G(F(t))-4G(t)\right]f(t)dt+G({\bf v})(1-2G({\bf v})).
\end{array} }
\end{equation}

Next, we observe that $H'(t)=f(F(t))F'(t)-f(t)$, and so the formula (\ref{eq3})  can be written as, after a change of variables in the second integral $t=F(t')$,

$$\begin{array}{c}
{\cal P}=\underset{\bf uv}{\int}G(F(t))\left[3-4H(t)\right]f(t)dt-\\ \\
-\underset{\bf uv}{\int}G(t)\left[3-4H(t)\right]f(F(t))F'(t)dt+G({\bf v})(1-2G({\bf v})) \Leftrightarrow
\end{array}$$

$$\begin{array}{c}
{\cal P}=\underset{\bf uv}{\int}H(t)\left[3-4H(t)\right]f(t)dt-\\ \\
-\underset{\bf uv}{\int}G(t)\left[3-4H(t)\right]H'(t)dt+G({\bf v})(1-2G({\bf v})).
\end{array}$$

\n An integration by parts in the last integral gives

$$\underset{\bf uv}{\int}G(t)\left[3-4H(t)\right]H'(t)dt=(3H(t)-2H(t)^2)G(t)|_{\bf u}^{\bf v}-\underset{\bf uv}{\int}H(t)\left[3-2H(t)\right]f(t)dt$$

\n Since $H({\bf u})=G({\bf v})$ and $H({\bf v})=1-G({\bf v})$, we obtain (\ref{equation1}). The function $H(t)$ (which under our assumptions must be continuous)
we see that this function must take  the value of $1/2$ by the intermediate value theorem. Indeed, $H(t)+H(F(t))=1$ and then if for some point $H(t_0)<1/2$ then
$H(F(t_0))>1/2$ which means somewhere, $H$ it must take the value $1/2$. Then we take that point to be $\bf u$. So the first term in formula (\ref{equation1}) is equal to $-1/2$, which shows  (\ref{equation2}) is true.

To obtain (\ref{equation3ind}), we apply (\ref{equation1}) for $u$, then for  $v=F(u)$ instead of $u$, and add the two identities up:
$$2{\cal P}=-(Q(G(v))+Q(1-G(v)))+6\underset{ \bf \partial R}{\int}H(t)\left[1-H(t)\right]f(t)dt,$$
where $Q(x)=x^2(3-2x)$. Since $Q(x)+Q(1-x)=1$, we get (\ref{equation3ind}).\eproof

\section{A few consequences}\label{consequences}
Let us continue with a result which seems to be known (see \cite{schneider&Weil} and \cite{wagner&welzl}). We will say that a region $\bf R$ is symmetric with respect to $O$ if $r(\theta)=r(\theta+\pi)$ for
all $\theta\in [0,2\pi]$.

\begin{theorem}\label{secondmain} The probability discussed in Theorem~\ref{main} satisfies  ${\cal P}\le 1/4$.
Moreover, ${\cal P}= 1/4$ if and only if $\bf R$ is symmetric.
\end{theorem}

\n \proof The first claim follows from (\ref{equation2}). If the region is symmetric with respect to $O$ then $H(t)=G(F(t))-G(t)=1/2$ for all $t\in {\bf uv}$.
By (\ref{equation2}), we have ${\cal P}=\frac{1}{4}$. Conversely, if ${\cal P}=\frac{1}{4}$ then $\underset{{\bf uv}}{\int}\left[\frac{1}{2}-H(t)\right]^2f(t)dt=0$ which attracts
$H(t)=1/2$ almost everywhere, but since we are dealing with at least continuous functions, we must have $H(t)=1/2$ for all $t$ on the boundary.  Hence $r(\theta+\pi)^2=r(\theta)^2$ which
implies that  $\bf R$ is symmetric with respect to $O$. $\Box$

\vspace{0.1in}

Formula (\ref{equation1}) is consistent with the intuitive fact stated next.

\begin{corollary}\label{convergencetozero} Suppose the point O approaches a point $\bf u$ on the boundary of $\bf R$ such that the region $\bf R$ is on one side of
a line through $\bf u$. Then  $$\lim_{O\to {\partial {\bf R} }}{{\cal P}_O}=0.$$
\end{corollary}
\n \proof If $O$ gets close to $\bf u$ on the boundary, there exists a line containing $O$ such that it cuts the region ${\bf R}$ into two regions, one of which has an area that goes to $0$
as $O$ gets closer to $u$. We can accomplish this by taking a parallel through $O$ to the line in the hypothesis of the corollary since the region $\bf R$ is bounded.
We take that line as ${\bf uv}$. Then by (\ref{equation1}), we have ${\cal P}_O\le \frac{3}{2}G({\bf v})\to 0$.$\Box$

We observe that this limiting behavior is not satisfied for every point on the boundary $\bf R$ (see Section~\ref{sectionOnLimacon}).

In \cite{Kleitman},  $\rm Hal\acute{a}sz$ and  Kleitman studied the maximum of $\cal P$ as a function of $O$ for a triangle. They show that this maximum is attainted for the center of mass
of the triangle. So, it is interesting to generalize this problem for arbitrary regions:

\vspace{0.1in}

{\it Given a fixed region $\bf R$, where is the position of $O$ to have a maximum of $\cal P_O$?}
\vspace{0.1in}

Theorem~\ref{secondmain} provides an answer to this question if $\bf R$ is symmetric but what happens in general?
The formula (\ref{equation2}) suggests that $\cal P$ is maximum where the function $t\to (1/2-H(t))$ has the smallest $L^2$ norm with respect to
this measure $f(x)dx$ on the boundary. On the other hand this measure depends on the point $O$  too, and as $O$ moves toward a point at which  $H(t)$
is closer to $1/2$ the density $f$ on $\bf uv$ gets bigger. So, the trade off is not clear and the problem can be decided only if a more
precise expression for ${\cal P}$ is found. However, the formula (\ref{equation1}) provides a better upper bound in the case of non-symmetric domains and this result seems to be new.

\begin{theorem}\label{upperboundnonsymm} The probability discussed in Theorem~\ref{main} satisfies
\begin{equation}\label{inequalitynonsymm}
\frac{1}{4}-(1+4h)\left(\frac{1}{2}-h\right)^2 \le  {\cal P}\le \frac{1}{4}-2\left(\frac{1}{2}-h\right)^3
\end{equation}
\n where $h={\bf \min}\{ H(t)| t \in {\bf \partial R}\}$.
\end{theorem}

\n \proof For the right-hand side of (\ref{inequalitynonsymm}), we observe that $H(t)(1-H(t))\le \frac{1}{4}$ for every $t$. Hence, using (\ref{equation1}), we get

$${\cal P}\le \frac{6}{4}G({\bf v})-G({\bf v})^2[3-2G({\bf v})]=\frac{x(3-6x+4x^2)}{2},$$

\n where $x=G({\bf v})$. Since $\bf u$ is arbitrary, this inequality is true for the smallest value of the function
$$W(x)=\frac{x(3-6x+4x^2)}{2}=\frac{1}{4}-\frac{(1-2x)^3}{4}=\frac{1}{4}-2\left(\frac{1}{2}-x\right)^3$$ over the interval
$[h,1/2]$ (we know that $x$ can be 1/2). The derivative of this function is $W'(x)=\frac{3(1-2x)^2}{2}\ge 0$ and so $W$ is strictly increasing. This gives a minimum of $W$
at $h$. Therefore, the right-hand side of (\ref{inequalitynonsymm}) follows from these simple considerations.

For the left-hand side of (\ref{inequalitynonsymm}), we use (\ref{equation11}) and observe that

$${\cal P}\ge  \frac{1}{4}-\frac{(1-2x)^3}{4} -6\left(\frac{1}{2}-h\right)^2x:=Z(x).$$

\n Since this inequality is true for every $u$, ($x=G(\bf v)$), we can take the maximum of the function $Z$.
Looking at the derivative of $Z$, we get $Z(x)=6(x-h)(x+h-1)$. Since $x\in [h,1-h]$ we see that $Z(x)\le 0$ and so
$Z$ is a decreasing function. So the best inequality we can obtain is for $x=h$, and then we obtain the left-hand side of (\ref{inequalitynonsymm}).
\eproof
\vspace{0.1in}

Another observation here is that $H$ has a minimum at a point $t_0$ where $H'(t_0)=0$. If we use polar representation, this means $r(t_0+\pi)=r(t_0)$.
In other words, $O$ divides the segment $t_0F(t_0)$ into equal parts.

For $\bf R$ a triangle and if $O$ is the center of mass, Theorem~\ref{upperboundnonsymm} shows that $\frac{176}{729}< {\cal P}\le \frac{182}{729}$
which is a decent approximation for the exact value of ${\cal P}$ computed in Section~5.

In the remaining sections we will show how our formula works in various situations where the calculations can be done and the result can be expressed in
a relatively simple way. 

\vspace{0.2in}
\section{A family of Lima\c cons between a circle and a  Cardioid}\label{sectionOnLimacon}
Perhaps the simplest curve that fits as the boundary of a region in our framework and for which all the involved calculations are quite elementary is the case of a Lima\c con.
We will start with the equations in polar coordinates $r=\frac{\sqrt{2}}{\sqrt{(2a^2+1)\pi}}(a+\cos t)$, $t\in [0,2\pi]$. This is, in fact, a family of  Lima\c cons over the parameter $a\in (1,\infty)$.
In Figure~4, we have included these curves for $a\in \{10,3,2,1.5,1\}$. Let us remark that $a=1$ is in fact a Cardioid and the point $O$ is on the boundary in this case.
At the other extreme, if $a\to \infty$ we get a circle centered at the origin of radius $\frac{1}{\sqrt{\pi}}$ (area 1).
\[
\underset{Figure\ 4}{\epsfig{file=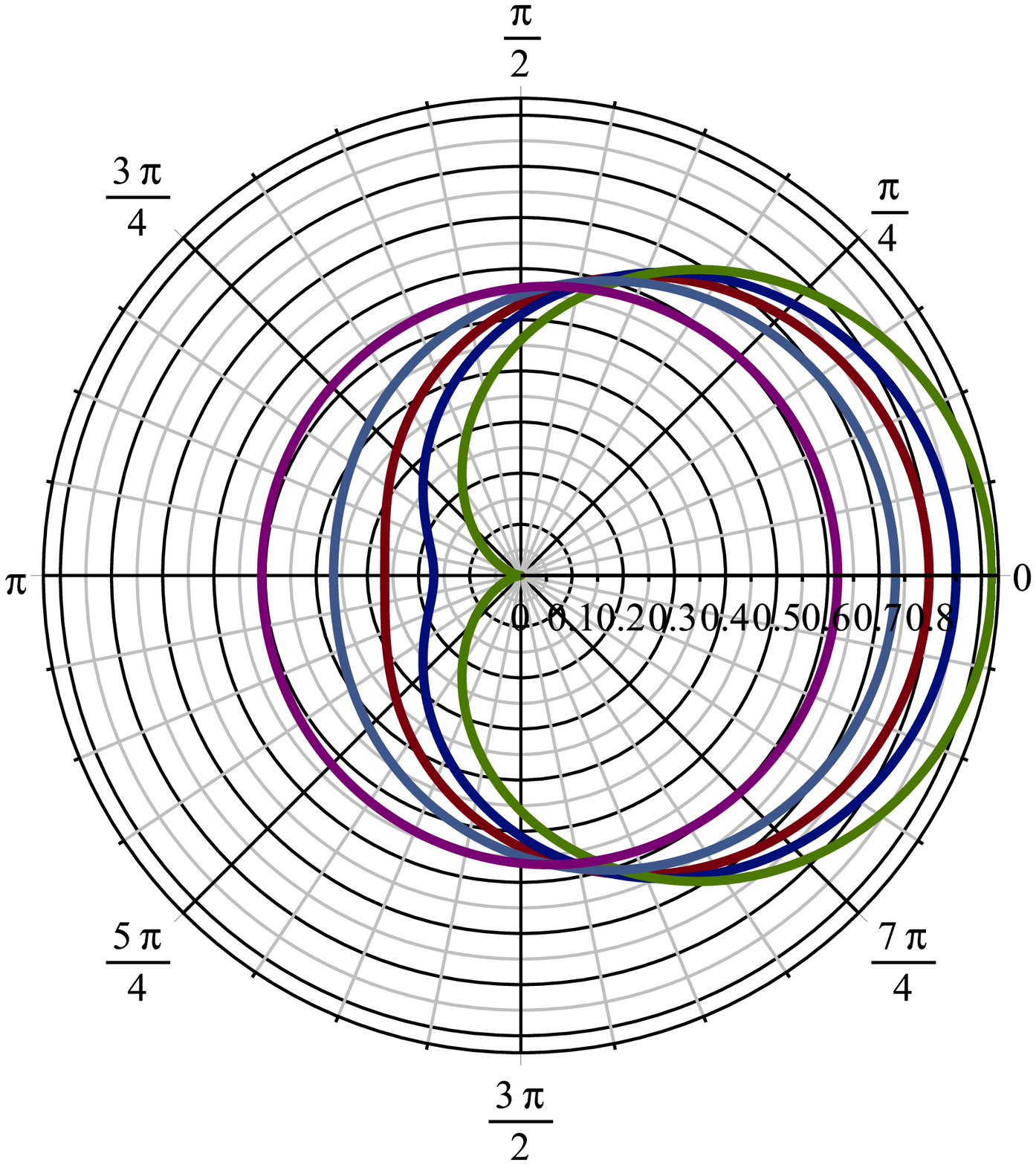,height=3.3in,width=3in}}
\]
We observe first that the area of each such curve is $1$:

$$Area=\frac{1}{2}\int_0^{2\pi} r(t)^2 dt=\frac{1}{(2a^2+1)\pi} \int_0^{2\pi} (a^2+2a\cos t+\frac{1+\cos 2t}{2}) dt=1. $$

So, the distribution on the boundary can be taken to be $\frac{1}{2}r(t)^2dt$. Then the function $H(t)=\int_t^{t+\pi} \frac{1}{2}r(t)^2dt$ or
$$H(t)=\frac{1}{(2a^2+1)\pi} \int_t^{t+\pi} (a^2+2a\cos t+\frac{1+\cos 2t}{2}) dt=\frac{a^2\pi-4a\sin t+\pi/2}{(2a^2+1)\pi} \Rightarrow$$

$$\frac{1}{2}-H(t)=\frac{4a\sin t}{(2a^2+1)\pi}, t\in[0,\pi].$$

Then, sine $H(0)=1/2$, formula (\ref{equation2}) gives
$${\cal P}=\frac{1}{4}-6\frac{16a^2}{(2a^2+1)^3\pi^3}\int_0^{\pi} (a^2+2a\cos t+\frac{1+\cos 2t}{2})\sin^2 tdt.$$

Since $(a^2+2a\cos t+\frac{1+\cos 2t}{2})\sin^2 t=a^2\frac{1-\cos 2t}{2}+\frac{2a}{3}\frac{d}{dt}(\sin ^3 t)+\frac{1-\cos^2 2t}{4}$ the above integration
leads to

$${\cal P}_a=\frac{1}{4}-6\frac{16a^2}{(2a^2+1)^3\pi^3}(\frac{a^2\pi}{2}+\frac{\pi}{8})=\frac{1}{4}-\frac{12a^2(4a^2+1)}{(2a^2+1)^3\pi^2}.$$

Interestingly enough as $a\to 1$ we get ${\cal P}_1=\frac{1}{4}-\frac{20}{9\pi^2}\approx 0.0248418142$. So, it is possible to have the probability ${\cal P}_O$
strictly positive for a point $O$ on the boundary of $\bf R$. We observe that
$$h=min H(t)=\frac{a^2\pi-4a+\pi/2}{(2a^2+1)\pi}=\frac{1}{2}-\frac{4a}{(2a^2+1)\pi}.$$

\n This checks out the inequality (\ref{inequalitynonsymm}) which reduces in this case to the obvious inequality  $\pi\ge \frac{32}{3(4a^2+1)}$ for $a>1$.

\section{The case when $\bf R$ is an equilateral triangle}\label{sectionEqTri}

\[
\underset{Figure\ 5}{\epsfig{file=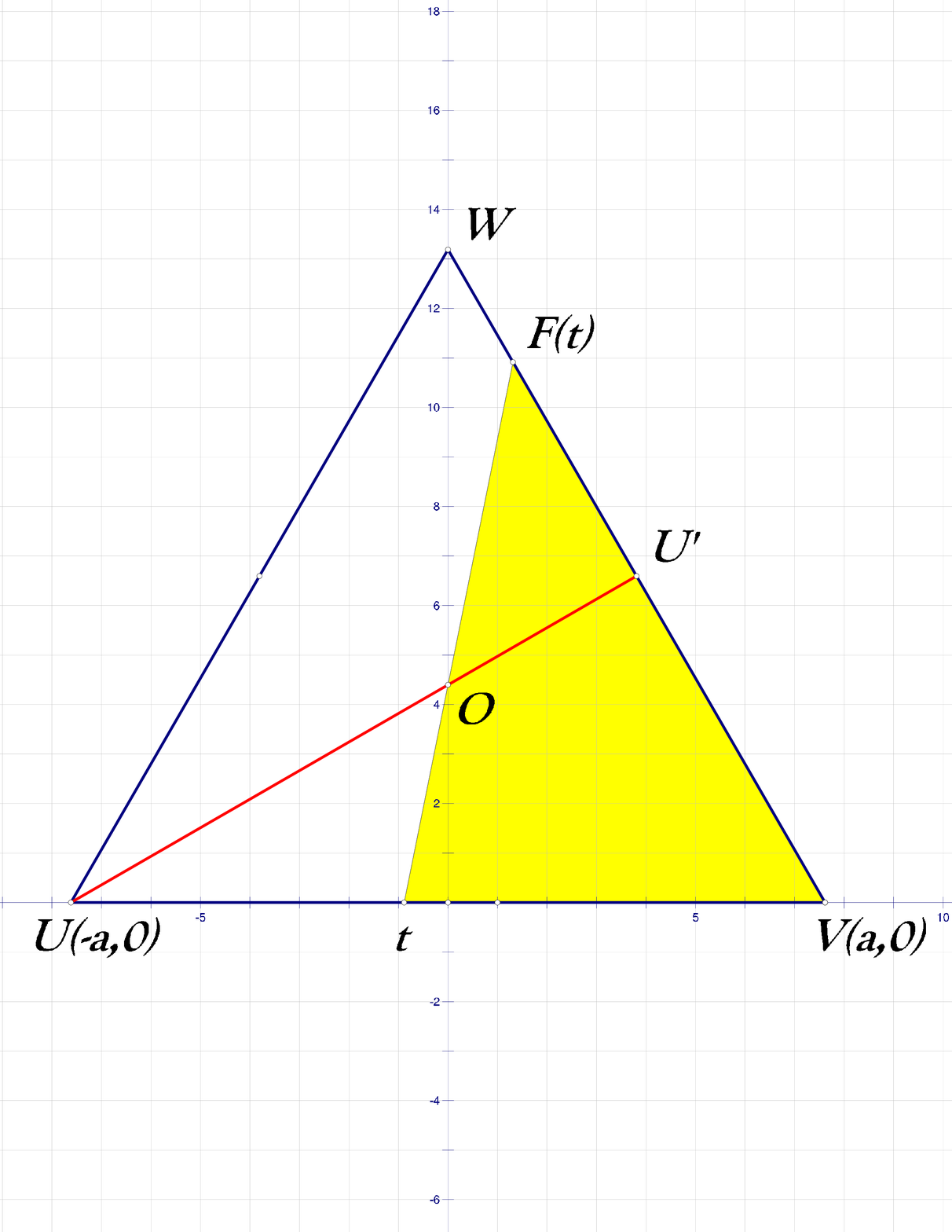,height=3.3in,width=3in}}
\]
Let us use formula (\ref{equation2}) to compute this probability for an equilateral triangle. We will do the calculations for the triangle positioned
as in Figure~5, with vertices $U(-a,0)$, $V(a,0)$ and $W(0,a\sqrt{3})$. The area of this triangle is $a^2\sqrt{3}$ so we will take $a=3^{-1/4}$ in order for the area to be equal to 1.
Since $UU'$ is a line of symmetry, we can use formula (\ref{equation2}). The distributions on all sides are all equal:  $f\equiv \frac{a\sqrt{3}}{6}=3^{-3/4}/2$.
To compute $F$ let us consider $t\in [-a,0]$ and observe that the line $tO$ has equation $X(\frac{a\sqrt{3}}{3}-0)-Y(0-t)-\frac{at\sqrt{3}}{3}=0$ and the line $VW$ has equation $X/a+Y/(a\sqrt{3})=1$.
Solving this system gives us $F(t)=(\frac{-2at}{a-3t},\frac{a(a-t)\sqrt{3}}{a-3t})$. This implies that $F(t)V=\frac{2a(a-t)}{a-3t}$. Hence,
$$H(t)=\frac{a(a-t)^2}{a-3t}\sin 60^{\circ}=\frac{a(a-t)^2\sqrt{3}}{2(a-3t)}$$ for $t\in [-a,0]$. Clearly, for $t\in [0,a]$ we have $H(t)=1-\frac{a(a+t)^2\sqrt{3}}{2(a+3t)}$
and so the function $H(t)(1-H(t))$ is symmetric with respect to each vertex and midpoint on the boundary $\partial R$.
As a result the calculation in (\ref{equation2}) reduces to

$${\cal P}=\frac{1}{4}-6\int_{U}^{U'}(\frac{1}{2}-H(t))^2f(t)dt=\frac{1}{4}-6\cdot 3\frac{a\sqrt{3}}{6}\int_{-a}^0(\frac{1}{2}-H(t))^2dt,$$

\n Changing the variable in the last integration $s=(a-t)/a$, $ds=-dt/a$, $s\in [1,2]$, we obtain

$${\cal P}=\frac{1}{4}-3a^2\sqrt{3}\int_{1}^2\left[\frac{1}{2}-\frac{s^2}{2(3s-2)}\right]^2ds=\frac{1}{4}-\frac{3}{4}\int_{1}^2\frac{(2-s)^2(s-1)^2}{(3s-2)^2}ds.$$

\n Changing the variable again, $3s-2=x$, the partial fraction decomposition comes easily:

$$\begin{array}{c}{\cal P}=\frac{1}{4}-\frac{1}{324}\int_1^4\frac{(4-x)^2(x-1)^2}{x^2}dx=\frac{1}{4}-\frac{1}{324}\int_1^4(x^2-10x+33-\frac{40}{x}+\frac{16}{x^2})dx=\\ \\
=\frac{1}{4}-\frac{1}{324}(\frac{x^3}{3}-5x^2-\frac{16}{x}+40\ln x)|_1^4=\frac{2(3+10\ln 2)}{81}\approx 0.2452215261 .
\end{array}$$

\n This is the same answer obtained in \cite{prekopa}.

\section{Regular Polygon with $2m+1$ vertices}\label{regularpolygon}
\[
\underset{Figure\ 6}{\epsfig{file=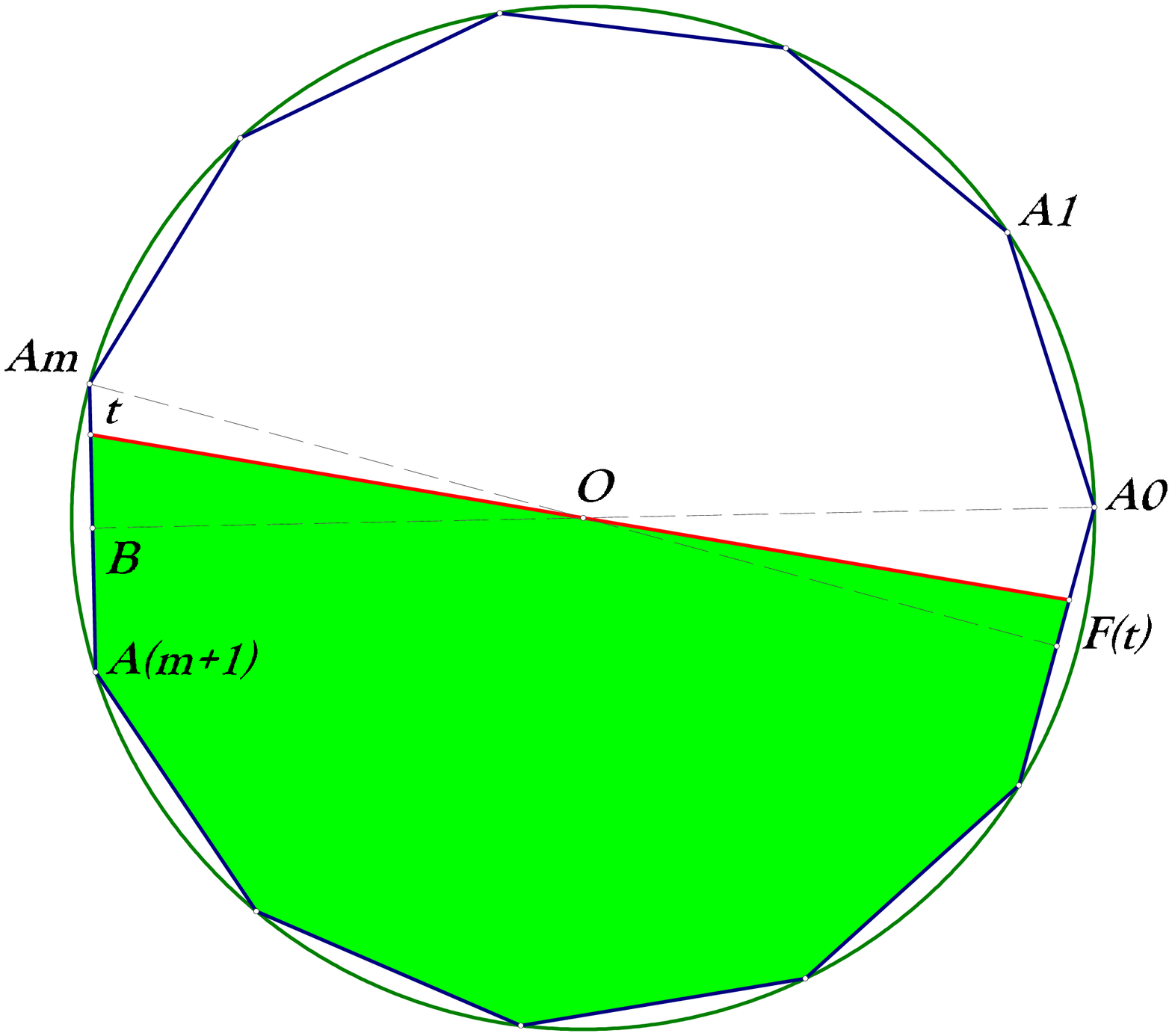,height=2.5in,width=3in}}
\]
Let us generalize the previous situation to a regular polygon with $2m+1$ vertices ($m\ge 1$). We will refer to Figure~5, where the regular polygon is centered at the origin
and it has vertices $A_k=Re^{i\theta}$ where $\theta_k=\frac{2k\pi}{2m+1}$, $k=0,1,...,2m$. Here $R$ is chosen so that the area of the polygon is equal to 1, i.e.,
$R=\sqrt{\frac{2}{(2m+1)\sin \frac{2\pi}{2m+1}}}$. The side lengths are $\ell=2R\sin (\theta_1/2)$. We will use that parametrization with respect to the arc-length along the boundary of
the polygon. As a result on each side $f(t)=\frac{1}{2} R\cos ((\theta_1/2))$. Let us introduce an usual notation here $a_p=R\cos ((\theta_1/2))$. We choose $\bf u$  to be $A_m=(-a_p,\ell/2)$ and as a result
${\bf v}=F({\bf u})$ is the midpoint of $\overline{A_{2m}A_0}$. The coordinates of ${\bf t}$ are $(-a_p,t)$, $t\in [0,\ell/2]$ and of $A_{2m}=(R\cos \theta_1,-R\sin \theta_1)$.
The equation of $\overline{\bf tF(t)}$ is $Y=-(t/a_p)X$ and the equation of $\overline{A_{2m}A_0}$ is $X(\sin \theta_1)-Y(1-\cos \theta_1)-R\sin \theta_1=0$, and so the intersection $\bf F(t)$
has coordinates $x_t=\frac{a_pR\sin \theta_1}{a_p\sin \theta_1+t(1-\cos \theta_1)}$ and $y_t=-\frac{tR\sin \theta_1}{a_p\sin \theta_1+t(1-\cos \theta_1)}$. Then,
$F(t)A_0=\frac{2tR\sin (\theta_1/2)}{ a_p\sin \theta_1+t(1-\cos \theta_1)}$. Therefore,

$$H(t)=\frac{1}{2}+Area({\bf t}BO)-Area(OA_0{\bf F(t)})=\frac{1}{2}(1+ta_p-a_pF(t)A_0)\Rightarrow$$

$$H(t)-\frac{1}{2}=\frac{a_p}{2}(t-\frac{2tR\sin (\theta_1/2)}{ a_p\sin \theta_1+t(1-\cos \theta_1)}).$$

Hence, we get
$${\cal P}=\frac{1}{4}-6(2m+1)\frac{a_p^3}{8}\int_0^{\ell/2}(t-\frac{2tR\sin (\theta_1/2)}{ a_p\sin \theta_1+t(1-\cos \theta_1)})^2dt.$$

Let us change the variable, $t=(\ell/2)s$, $s\in [0,1]$ and also we will denote $\theta=\theta_1/2$. This implies
$${\cal P}=\frac{1}{4}-3(2m+1)\frac{R^3\cos^3 \theta R^3\sin^3 \theta}{4}\int_0^1(s-\frac{2s\sin\theta}{ \cos \theta \sin 2\theta+s\sin \theta (1-\cos 2\theta)})^2ds,$$

\n or

$${\cal P}=\frac{1}{4}-\frac{3}{4(2m+1)^2}\int_0^1(s-\frac{s}{ \cos^2 \theta +s \sin^2 \theta})^2ds.$$

\n This integration is rather difficult one, but with some help from Maple one obtains

\begin{equation}\label{regpoly}
{\cal P}_m=\frac{1}{4}-\frac{1+9a-9a^2-a^3+6a(1+a)\ln a}{4(2m+1)^2(1-a)^3},
\end{equation}

\n where $a=\cos^2 \theta=\cos ^2\frac{\pi}{2m+1}$.

\n For a regular pentagon, $m=2$ and $a=\frac{3+\sqrt{5}}{8}$, one gets a pretty complicated expression for the probability but it is an exact one:

$${\cal P}_5=\frac{1}{625}\left[(990+438\sqrt{5})\ln (\sqrt{5}-1)-30(2+3\sqrt{5})\right]\approx  0.24982224 .$$

\n It is not difficult to check that if $m\to \infty$, the function in (\ref{regpoly}) converges indeed to $1/4$.

\section{An arbitrary triangle and an arbitrary point in its interior}
\[
\underset{Figure\ 7}{\epsfig{file=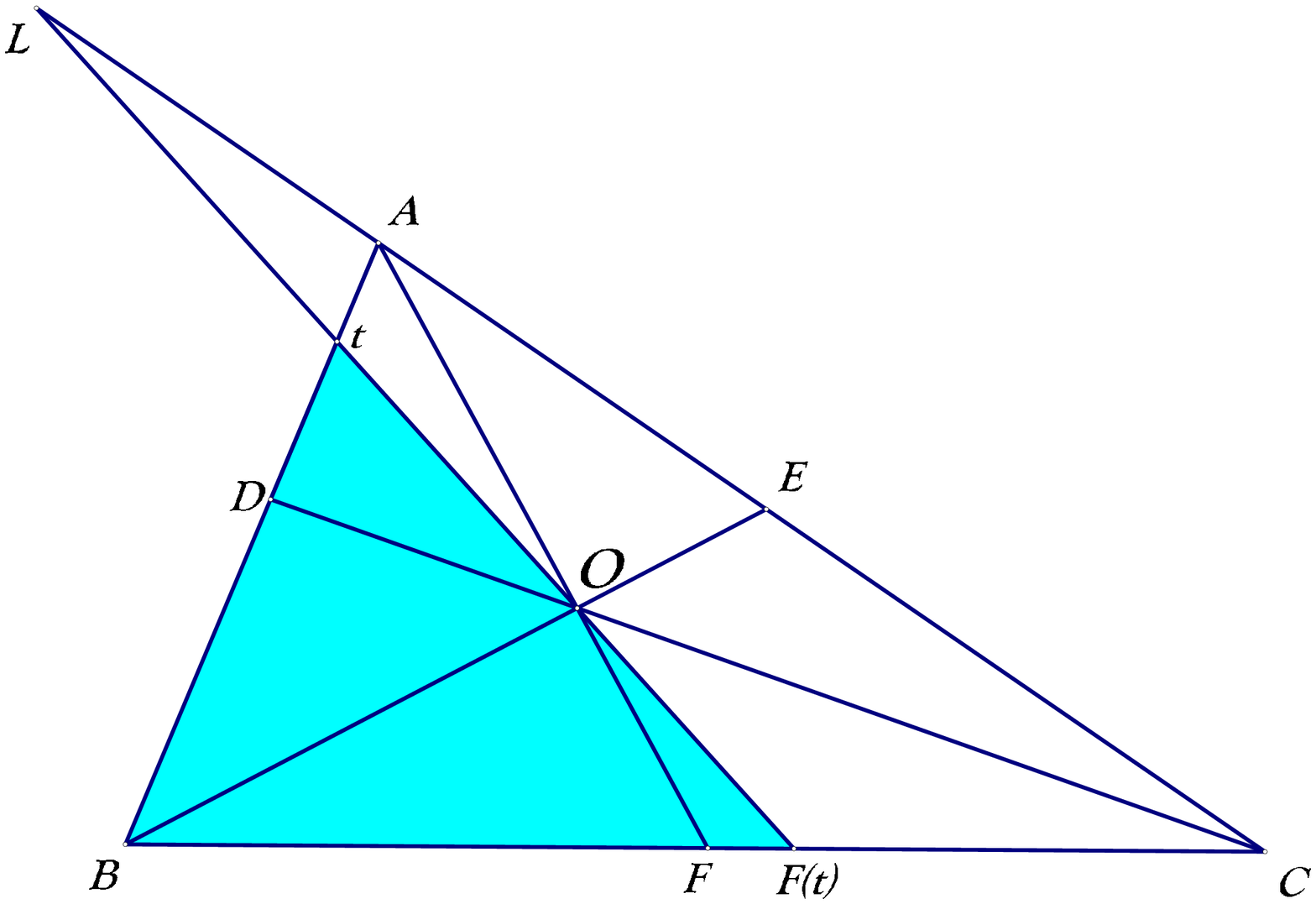,height=2.5in,width=3in}}
\]
Let us start with an arbitrary triangle $ABC$ of area equal to 1 (as in Figure~7), and an arbitrary point $O$ in the interior of the triangle.
We let $D=F(C)$, $E=F(B)$ and $F(A)=F$. For  $\overline{\bf uv}$ we take  $\overline{AF}$ and then formula (\ref{equation1})
implies:

$$
{\cal P}=-G({\bf v})^2[3-2G({\bf v})]+6\underset{\ \ {\bf \overarc{\bf uv}}}{\int}H(t)\left[1-H(t)\right]f(t)dt.
$$
 In this situation $G({\bf v})$ is just the area of the triangle $AFC$.
So, we will integrate over the segments $\overline{AD}$, $\overline{DB}$ and $\overline{BF}$. Let us denote the distances from $O$ to the sides
$\overline{BC}$, $\overline{CA}$, and $\overline{AB}$ by $d_A$, $d_B$ and $d_C$ respectively. Also, we introduce the standard notation for the
altitudes of the triangle $h_A$, $h_B$ and $h_C$. This gives us the barycentric coordinates of $O$: $\alpha=\frac{d_A}{h_A}$, $\beta=\frac{d_B}{h_B}$
and $\gamma=\frac{d_C}{h_C}$. We have

$$\alpha+\beta+\gamma=\frac{(BC\cdot d_A)/2}{(BC\cdot h_A)/2}+\frac{(AC\cdot d_B)/2}{(AC\cdot h_B)/2}+\frac{(AB\cdot d_C)/2}{(AB\cdot h_C)/2}=[OBC]+[OAC]+[OAB]=1,$$

\n where $[XYZ]$ denotes the area of the triangle $XYZ$.

The distribution $f$ along the sides is piecewise constant: $\frac{d_C}{2}$ on $\overline{AB}$,
$\frac{d_A}{2}$ on $\overline{BC}$, and $\frac{d_B}{2}$ on $\overline{AC}$.

First, we observe that

$$\frac{DA}{DB}=\frac{[ADO]}{[BDO]}=\frac{[ADC]}{[BDC]}=\frac{[AOC]}{[BOC]}=\frac{\beta}{\alpha},$$

\n and similarly $$\frac{EA}{EC}=\frac{\gamma}{\alpha},\ \ \ \text{and}\ \ \frac{FB}{FC}=\frac{\gamma}{\beta}.$$

Then,
$$G({\bf v})=[ABF]=\frac{[ABF]}{[ABC]}=\frac{BF}{BC}=\frac{\gamma}{\beta+\gamma}.$$

In order to calculate the function $H(t)$ for $t\in \overline{AD}$, we denote by $x$ the length of the segment $\overline{At}$.

\begin{lemma}\label{menelaus}
For $t\in \overline{AD}$, we have
\begin{equation}\label{identityTriangle}
\alpha \frac{tA}{tB}+\gamma \frac{F(t)C}{F(t)B}=\beta.
\end{equation}
\end{lemma}

\n \proof If $\overline{tF(t)}$ is parallel to $\overline{AC}$ then the two ratios in (\ref{identityTriangle}) are equal and equal to $OE/OB=\frac{d_B}{h_b-d_B}=\frac{\beta}{1-\beta}$, which is precisely (\ref{identityTriangle}). If $\overline{tF(t)}$ is not parallel to $\overline{AC}$, let $L$ be their intersection (see Figure 6). Since  $\frac{EA}{EC}=\frac{\gamma}{\alpha}$ we have $\frac{LE-LA}{LC-LE}=\frac{\gamma}{\alpha}$.
Solving for $LE$ from this last equality we obtain $LE=\frac{\alpha LA+\gamma LC}{\alpha+\gamma}$. By Menelaus' Theorem in the triangle $ABE$ and the transversal $L-t-O$ we get
$\frac{tA}{tB}\frac{OB}{OE}\frac{LE}{LA}=1$ and solving for $\frac{tA}{tB}=\frac{LA}{LE}\frac{\beta}{1-\beta}$. Similarly, Menelaus' Theorem in the triangle
$BOC$ and the transversal $L-O-F(t)$ gives $\frac{F(t)C}{F(t)B}=\frac{LC}{LE}\frac{\beta}{1-\beta}$.
Then, using the relation about $LE$ that we established earlier we get

$$\alpha \frac{tA}{tB}+\gamma \frac{F(t)C}{F(t)BC}=\frac{\alpha LA+\gamma LC}{LE}\frac{\beta}{1-\beta}= (\alpha+\gamma)\frac{\beta}{1-\beta}=\beta,$$

\n which proves (\ref{identityTriangle}). \eproof

Lemma~\ref{menelaus} allows us to calculate $H(t)$ as a function of $x$, since the area of the triangle
$tBF(t)$ can be written as
$$[tBF(t)]=\frac{ [tBF(t)]}{[ABC]}=\frac{tB}{AB}\frac{F(t)B}{BC}=(1-\frac{x}{AB})(1-\frac{F(t)C}{BC}).$$

From (\ref{identityTriangle}) we obtain

$$\frac{F(t)C}{BC}=\frac{\beta-\alpha \frac{x}{AB-x}}{\beta+\gamma -\alpha \frac{x}{AB-x}}=
\frac{\beta-(\alpha+\beta)\frac{x}{AB}}{\beta+\gamma-\frac{x}{AB}} \Rightarrow $$

$$1-\frac{F(t)C}{BC}=\frac{\gamma (1-\frac{x}{AB})}{\beta+\gamma-\frac{x}{AB}}. $$

This implies that  $H(t)=\frac{\gamma (1-\frac{x}{AB})^2}{\beta+\gamma-\frac{x}{AB}}$. In the integral over $\overline{AD}$ we make a substitution $x=AB s$, with $s\in [0,\frac{\beta}{\alpha+\beta}]$ and hence, the contribution of this integral is

$$\underset{\ \ {\bf \overline{\bf AD}}}{\int} H(t)\left[1-H(t)\right]f(t)dt=
AB (d_C/2)\int_0^{\frac{\beta}{\alpha+\beta}}\frac{\gamma (1-s)^2}{\beta+\gamma-s}-\frac{\gamma^2 (1-s)^4}{(\beta+\gamma-s)^2}ds.$$

Therefore, this reduces to

$$\underset{\ \ {\bf \overline{\bf AD}}}{\int} H(t)\left[1-H(t)\right]f(t)dt=\gamma ^2 \int_0^{\frac{\beta}{\alpha+\beta}}\frac{ (1-s)^2}{\beta+\gamma-s}-\frac{\gamma (1-s)^4}{(\beta+\gamma-s)^2}ds.$$

If we change the variable again, $\beta+\gamma-s=u$, we obtain

$$I(\alpha,\gamma)=\gamma^2\int_{\frac{\alpha \gamma}{1-\gamma}}^{1-\alpha}\frac{(\alpha +u)^2}{u}-\frac{\gamma (\alpha +u)^4}{u^2} du.$$

The integrals over $\overline{\bf DB}$ and $\overline{\bf BF}$ are written in a similar way, and so we get an integral formula which for computational purposes
is very easy to use even by hand:

\begin{equation}\label{myeqtri}
{\cal P}=6(I(\alpha,\gamma)+I(\beta,\gamma)+I(\beta,\alpha))-\frac{\gamma^2(3-3\alpha-2\gamma)}{(1-\alpha)^3}.
\end{equation}

\n For instance if $\alpha=\frac{1}{6}$, $\beta=\frac{1}{3}$ and $\gamma=\frac{1}{2}$ then ${\cal  P}=\frac{1}{27}+\frac{41\ln 5}{972}+\frac{17\ln 2}{243} \approx 0.1534$.
We used formula (\ref{myeqtri}) to plot this function in 3D over the points of an equilateral triangle in terms of the barycentric coordinates of the points inside the equilateral triangle:

\[
\underset{Figure\ 8}{\epsfig{file=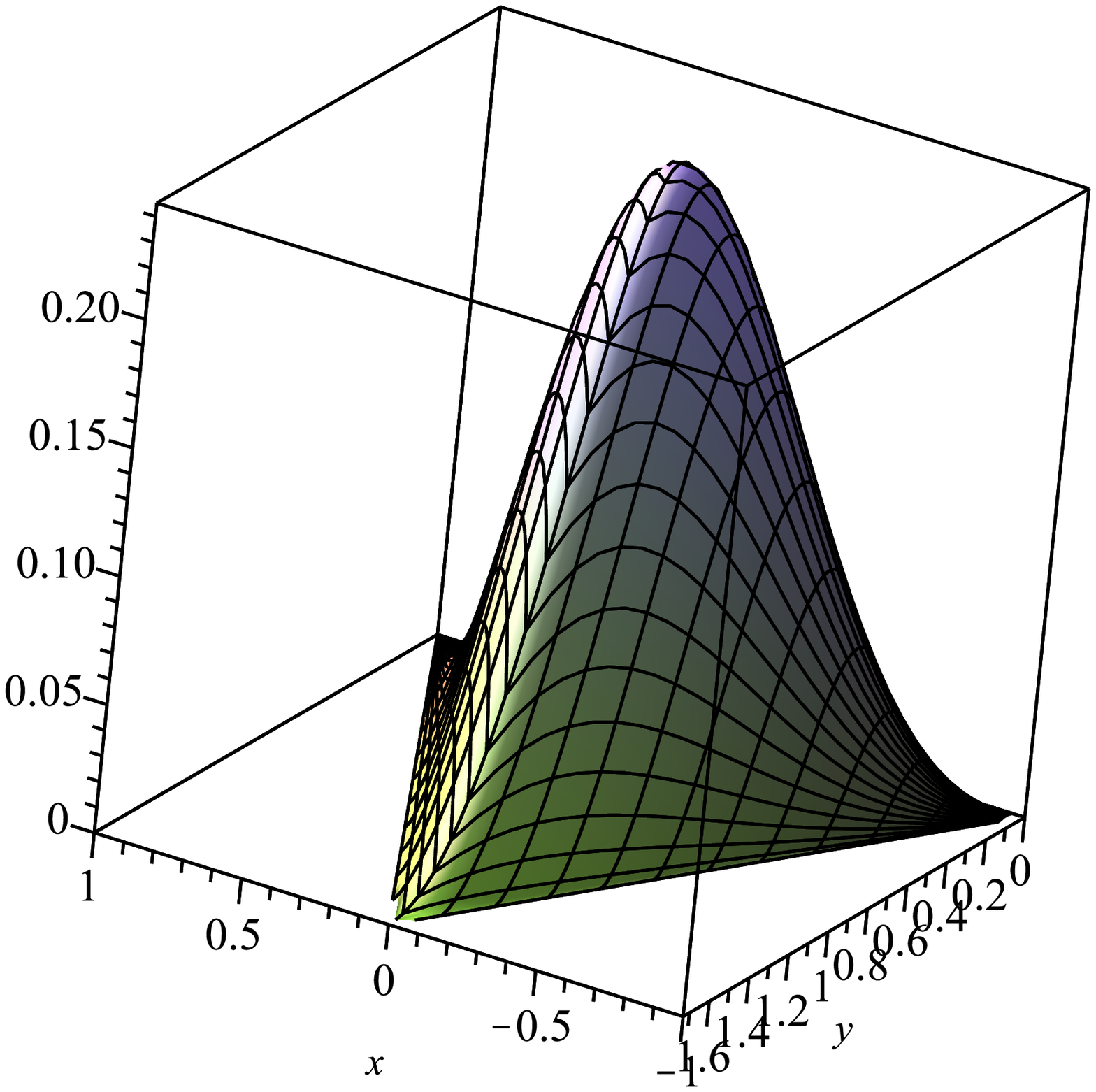,height=2.5in,width=3in}}
\]
Clearly this graph is consistent to the result that the maximum of the probability in question is attained for $\alpha=\beta=\gamma=\frac{1}{3}$.
It is interesting that one can now use (\ref{myeqtri}) to compute $\cal P$ for $O$ the incenter of triangle, the circumcenter, or the orthocenter.
For instance if $O=I$, and the sides of the triangle are $a$, $b$ and $c$, then taking $\alpha=\frac{a}{a+b+c}$, $\beta=\frac{b}{a+b+c}$, and $\gamma =\frac{c}{a+b+c}$.

\section{The case of a square and an arbitrary point}

 In what follows we will refer to Figure~9. Without loss of generality we
assume the square to be the unit square $S:=[0,1]^2$ and the point $O(u,v)$ with $u,v\in (0,1)$.  The distribution of the points $A'=s(A)$ on the sides is almost uniform.
Indeed, let us say the distribution of $A'$ on side $UV$ is given by a continuous function $f$. Then for $0\le a\le b\le 1$ we have

$$\int_a^b f(x)dx={\mathbb P}(A'\in [a,b])=\text{Area of triangle}\  \triangle OA'B'=\frac{1}{2}(b-a)\cdot v.$$

\n Differentiating with respect to $b$ gives $f(b)=\frac{v}{2}$. We denote this distribution corresponding to side $UV$ by $f_{UV}$ and similarly $f_{VW}(y)=\frac{1-u}{2}$, $f_{WZ}(x)=\frac{1-v}{2}$, and finally  $f_{ZU}(y)=\frac{u}{2}$. As we have noticed before, we have

$$\int_{\partial S}f(s)ds=\frac{v}{2}+\frac{1-u}{2}+\frac{1-v}{2}+\frac{u}{2}=1.$$

\[
\underset{Figure\ 9}{\epsfig{file=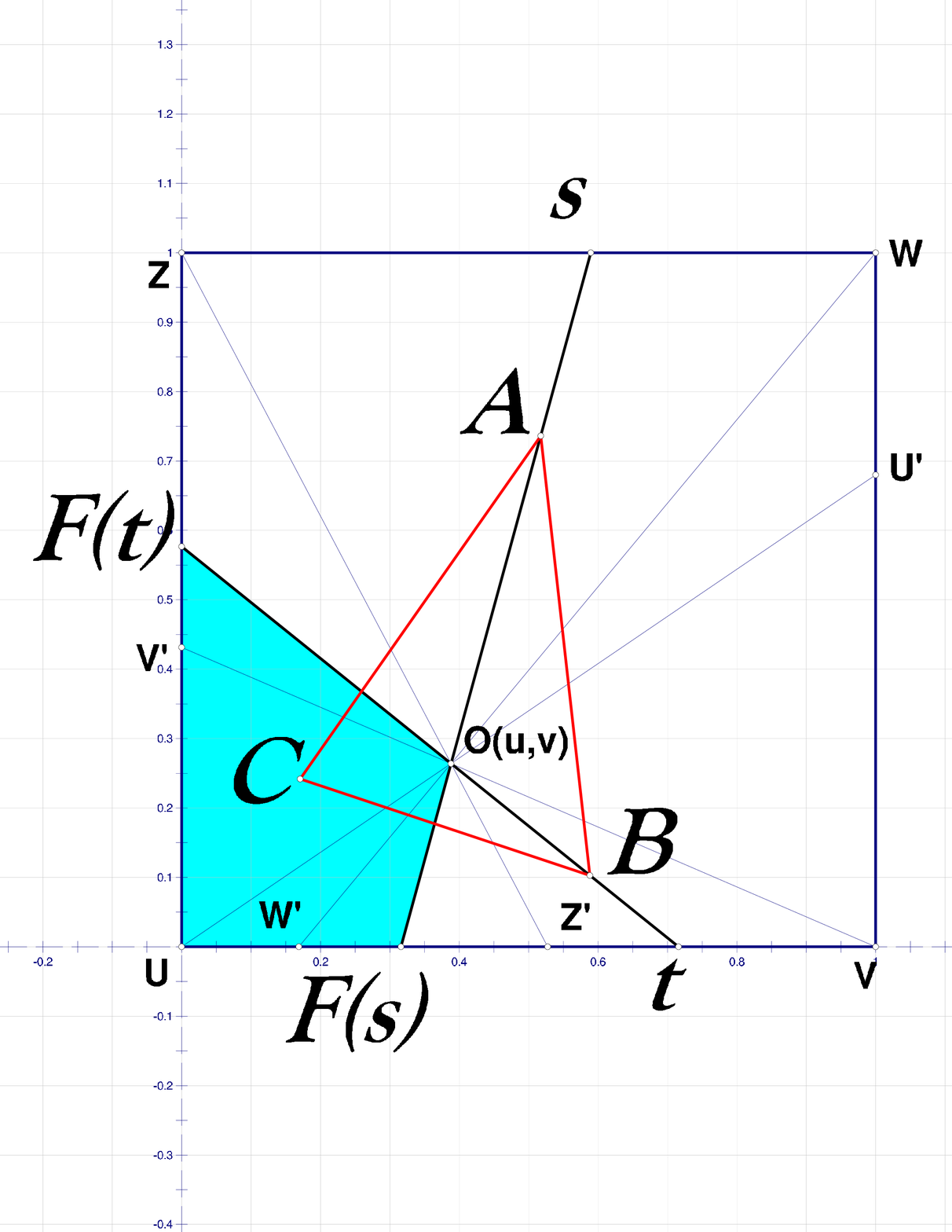,height=3.3in,width=3in}}
\]

Due to the symmetries of the square, we may assume, without loss of generality that $0<v< u<\frac{1}{2}$ (the limiting cases $u=v$ and $u=1/2$, can be obtained by a continuity argument).
We let $U'=F(U)$, $V'=F(V)$, $W'=F(W)$ and $Z'=F(Z)$.   These points, together with the vertices of the square, divide the boundary $\partial S$ into eight parts.
The coordinates of these points are easy to compute: $U'=(1,\frac{v}{u})$, $V'=(0,\frac{v}{1-u})$, $W'=(\frac{u-v}{1-v},0)$ and $Z'=(\frac{u}{1-v},0)$.
So, we have a piecewise definition for $F$ on $8$ intervals:

$$\small F(t)=\begin{cases} (1,\frac{v(1-x)}{u-x})\ \ \text{if}\ t=(x,0),x\in [0,\frac{u-v}{1-v}],\\ \\
(\frac{u-(1-v)x}{v},1)\  \ \text{if}\ t=(x,0),x\in [\frac{u-v}{1-v},\frac{u}{1-v}],\\ \\
(0,\frac{vx}{x-u})\  \ \text{if}\ t=(x,0),x\in [\frac{u}{1-v},1],\\ \\
(0,\frac{v-uy}{1-u}) \  \ \text{if}\ t=(1,y),y\in [0,\frac{v}{u}],
\end{cases}, \ \
\small F(t)=\begin{cases} (\frac{uy-v}{y-v},)\ \ \text{if}\ t=(y,0),y\in [\frac{v}{u},1],\\ \\
(\frac{u-vx}{1-v},0)\  \ \text{if}\ t=(x,1),x\in [0,1],\\ \\
(0,\frac{uy}{y-v})\  \ \text{if}\ t=(0,y),x\in [\frac{v}{1-u},1],\\ \\
(0,\frac{v-(1-u)y}{u}) \  \ \text{if}\ t=(0,y),y\in [0,\frac{v}{1-u}].
\end{cases}
$$

\n For the function $H$ we just need to compute some areas:
$$
\small H(t)=\begin{cases} \frac{v(1-x)^2}{2(u-x)}\ \ \text{if}\ t=(x,0),x\in [0,\frac{u-v}{1-v}],\\ \\
\frac{2vx+u-x}{2v}\  \ \text{if}\ t=(x,0),x\in [\frac{u-v}{1-v},\frac{u}{1-v}],\\ \\
\frac{vx^2}{2(x-u)}\  \ \text{if}\ t=(x,0),x\in [\frac{u}{1-v},1],\\ \\
\frac{v+y(1-2u)}{2(1-u)} \  \ \text{if}\ t=(1,y),y\in [0,\frac{v}{u}]
\end{cases},\ \
\small H(t)=\begin{cases} \frac{(1-u)y^2}{2(y-v)}\ \ \text{if}\ t=(y,0),y\in [\frac{v}{u},1],\\ \\
\frac{u+x-2vx}{2(1-v)}\  \ \text{if}\ t=(x,0),x\in [0,1],\\ \\
\frac{uy^2}{2(y-v)}\  \ \text{if}\ t=(x,0),x\in [\frac{v}{1-u},1],\\ \\
\frac{v-y+2uy}{2u} \  \ \text{if}\ t=(0,y),y\in [0,\frac{v}{1-u}].
\end{cases}
$$

\n Because there is no obvious point $\bf u$ for which $G(\bf v)=1/2$ we will use formula (\ref{equation1}).
Then the probability $\cal P$ is given by four integrals ${\cal P}=6(I_1+I_2+I_3+I_4)-\frac{v^2(3u-v)}{4u^3}$ where

 $$I_1=\int_0^{\frac{u-v}{1-v}}\frac{v^2(1-x)^2}{4(u-x)}(1-\frac{v(1-x)^2}{2(u-x)})dx,$$

 $$I_2=\int_{\frac{u-v}{1-v}}^{\frac{u}{1-v}}\frac{2vx+u-x}{4}(1-\frac{2vx+u-x}{2v})dx,$$

 $$I_3=\int_{\frac{u}{1-v}}^1\frac{v^2x^2}{4(x-u)}(1-\frac{vx^2}{2(x-u)})dx,\ \text{and}$$

 $$I_4=\int_0^{\frac{v}{u}}\frac{v+y(1-2u)}{4}(1-\frac{v+y(1-2u)}{2(1-u)})dy.$$

\n The answer in general is rather complicated and it is not symmetrical as expected:

$${\cal P}=Q_1(u,v)+Q_2(u,v)\ln \frac{1-u}{u}+Q_3(u,v)\ln \frac{1-v}{v}$$

\n where $Q_i$ are rational functions in $u$ and $v$.

\n For particular situations one can obtain pretty short expressions, for instance,
if $u=1/2$ and $v=1/4$ we obtain ${\cal P}=\frac{5}{48}+\frac{9\ln 3}{256}$.

In the case $u=v$, the answer is a little more manageable

\begin{equation}\label{squareu=v}
{\cal P}_u=\frac{1}{4}-\frac{(1-2u)(1-2u^2)(1+u-6u^3))}{4(1-u)}+3u^4(1-2u^2)\ln \frac{1-u}{u},\ \ u\in (0,\frac{1}{2}).
\end{equation}

If $u=v=\frac{1}{3}$ this reduces to $\frac{23}{162}+\frac{7\ln 2}{243}\approx 0.161942$.
We checked these experimentally and in $10^6$ trials they check usually with the third decimal.
As a curiosity, for $u=v=\frac{1}{1+e}$, the probability is just in terms of $e$: ${\cal P}=\frac{5e^5+6e^4+13e^3+7e^2-6e+1}{e(e+1)^6}$.
\section{A disk and its center, without a ``pizza" type slice }

The region we are interested in is included in Figure~12 and described by
\begin{equation}\label{circlewithoutapart}
\begin{array}{c}{\bf \cal R}_{r,\alpha}=\{(x,y)\in {\mathbb R^2} |x=t\cos \theta, y=t\sin \theta,\ \pi\ge |\theta|\ge \alpha,t\in [0,r] \},\\ \\ \ \ \alpha=\pi a, a\in [0,\frac{1}{2}],r>0,
\end{array}
\end{equation}

\n where the triangle $\triangle ABC$ contains the origin $O(0,0)$. We show that the correspoding probability is equal to

 $${\cal P}_a=\frac{(1+a)(1-2a)^2}{4(1-a)^3}.$$
 \[
\underset{Figure\ 12}{\epsfig{file=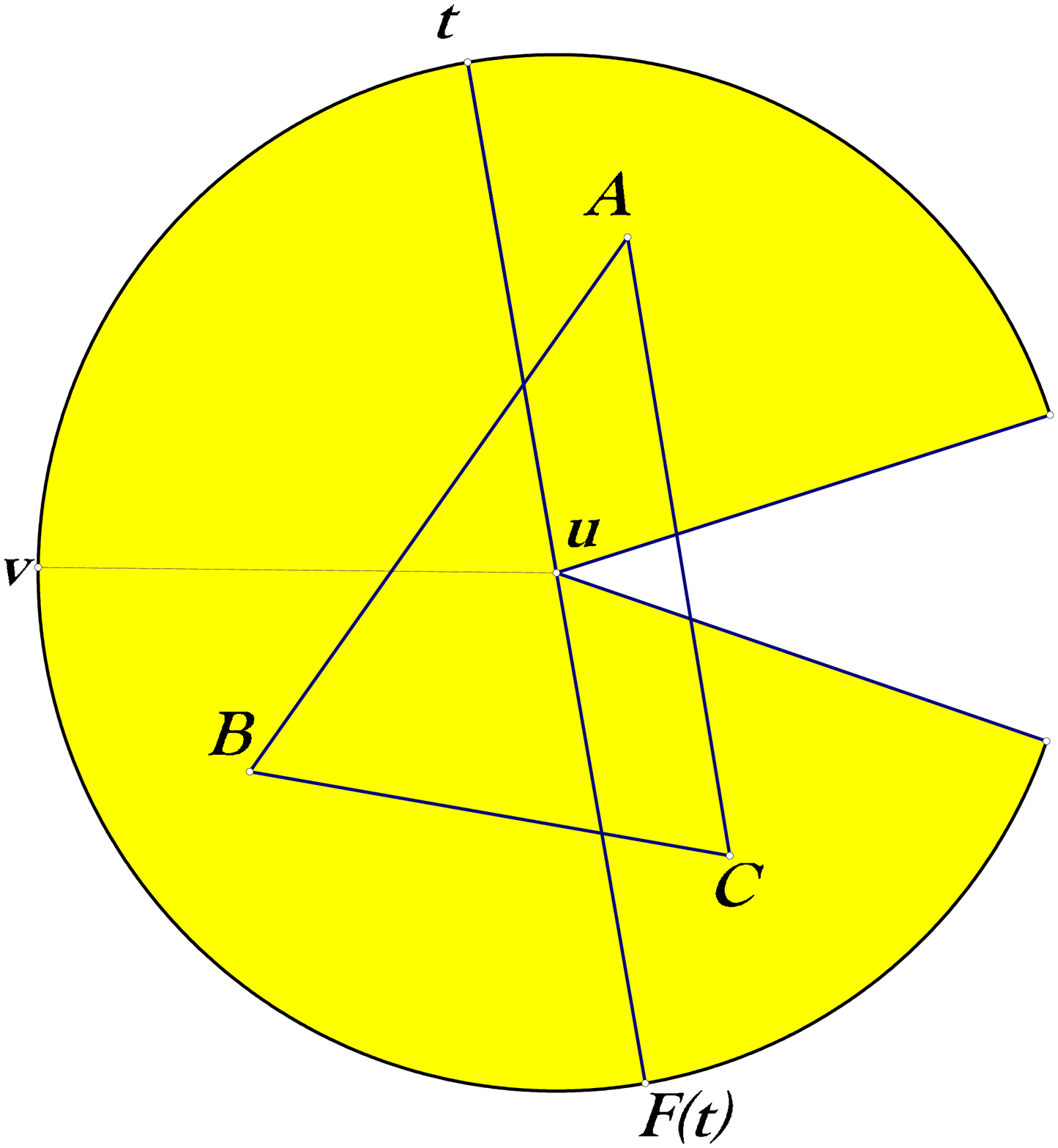,height=2.5in,width=2.5in}}
\]
 Without loss of generality we may assume that
the area of the region  ${\bf \cal R}$ is $1$, or $r=\frac{1}{\sqrt{\pi-\alpha}}$. We will use the formula (\ref{equation1}) for $u=O(0,0)$ and $v=(-r,0)$:

$${\cal P}=-G(v)^2(3-2G(v))+6\int_{uv}H(t)(1-H(t))f(t)dt,$$

\n where $u=(r,0)$. Hence, $G(v)=\frac{1}{2}$ and $H(t)=\frac{\pi}{2(\pi-\alpha)}$ for every $t=(r,\theta)$ with $\theta \in[\alpha,\pi -\alpha ]$.

The function $f$ is constant along the round boundary, and this constant is $\frac{1}{2(\pi-\alpha)}$ if we run the integration over the polar parameter $\theta\in [\alpha,2\pi -\alpha]$.
On the two segments, $f$ is clearly identically equal to $0$. We have $H(t)=\frac{1}{2}+\frac{\pi-\theta}{2(\pi-\alpha)}$ if $\theta\in [\pi-\alpha,\pi]$.
Then the probability is

$${\cal P}=-\frac{1}{2}+6\frac{1}{2(\pi-\alpha)}(\int_{\alpha}^{\pi-\alpha} \frac{\pi(\pi-2\alpha)}{4(\pi-\alpha)^2 }d\theta +
\int_{\pi-\alpha}^{\pi} (\frac{1}{4}-\frac{(\pi-\theta)^2}{4(\pi-\alpha)^2})d\theta ).$$

\n If we set $\alpha=\pi a$, $a\in [0,1/2]$, and changing the variable $\theta=\pi s$ we get

 $${\cal P}_a=-\frac{1}{2}+\frac{3(1-2a)^2}{4(1-a)^3}+\frac{3a}{4(1-a)}-\frac{a^3}{4(1-a)^3}=\boxed{\frac{(1+a)(1-2a)^2}{4(1-a)^3}}.\ \ \bsq$$

\section{Circle with a point off center}

\[
\underset{Figure\ 11}{ \epsfig{file=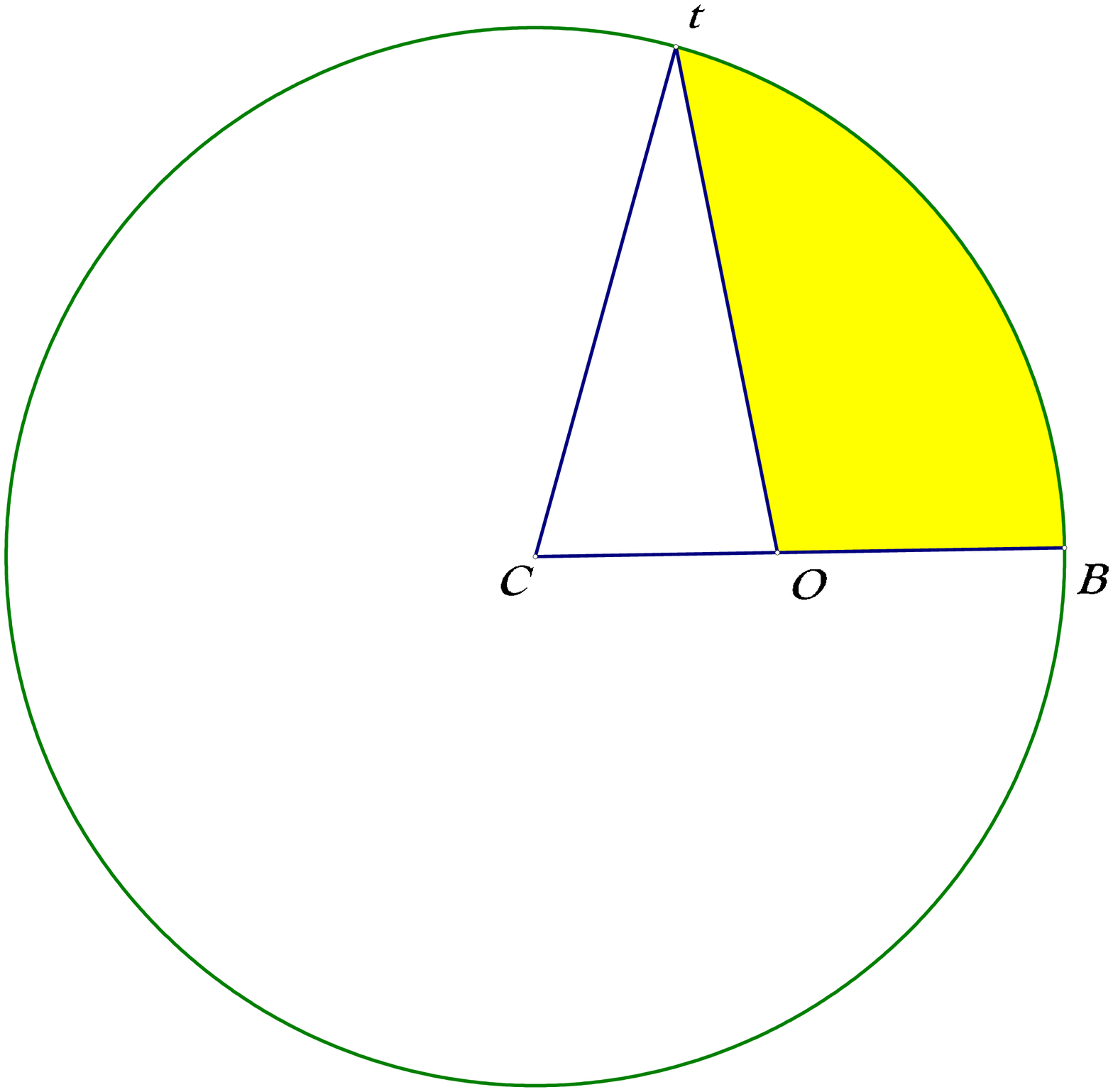,height=3in,width=3in}\ \epsfig{file=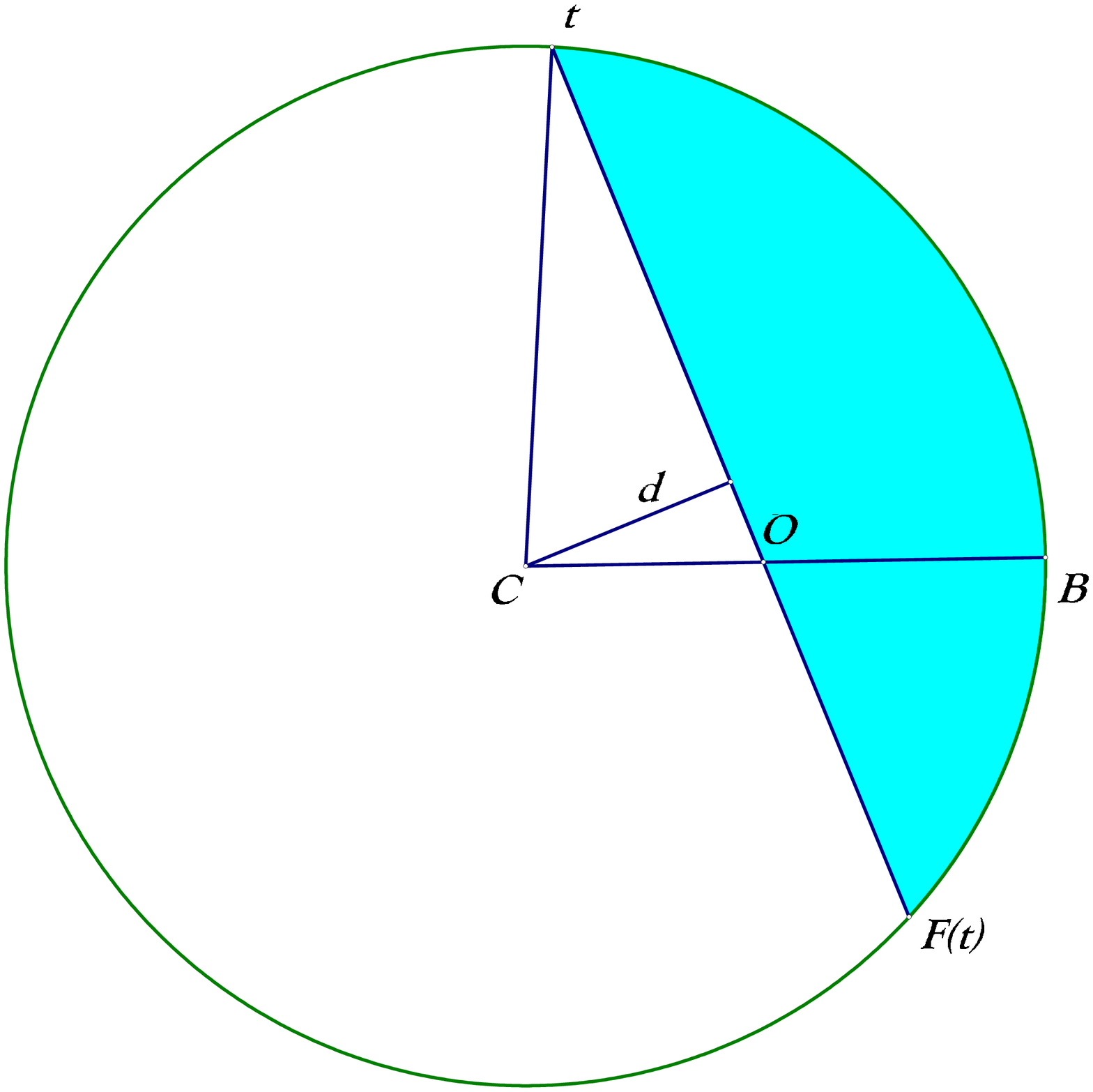,height=3in,width=3in}}
\]

Let us work with a circle centered at the origin of radius $R=\frac{1}{\sqrt{\pi}}$ (Figure~11). The point $O$ has coordinates $(rR,0)$ with $r\in [0,1)$.
Next, we will check that the function $f$ is given by $f(s)=R^2(1-r\cos s)/2$ in the sense that
$$\int_0^t f(s)ds=Area(O\overarc{Bt}), t\in [0,2\pi] $$
Indeed the above equality is the same as
$$\frac{1}{2}R^2(t-r\sin t)=Area(sector \ TCB)-Area(\triangle TCO)= Area(O\overarc{Bt}),$$
which is correct.  The function $H$ is given by a well known formula, $1-H(t)=R^2\arccos(\frac{d}{R})-d\sqrt{R^2-d^2}$, where $d$ is the distance from $C$ to the line
$tF(t)$: $d=\frac{R^2r\sin t}{tO}$. Since $tO=R\sqrt{1-2r\cos t+r^2}$ we obtain
$$H(t)=1-R^2\left[\arccos(\frac{r\sin t}{\sqrt{1-2r\cos t+r^2}})-\frac{r(1-r\cos t)\sin t }{1-2r\cos t+r^2}\right], t\in [0,\pi].$$

\n Using (\ref{equation2}), the probability in this case is given by

$${\cal P}_r=\frac{1}{4}-\frac{3}{\pi }\int_0^{\pi}\left[\frac{1}{2}-\frac{1}{\pi}\arccos(\frac{r\sin t}{\sqrt{1-2r\cos t+r^2}})+\frac{r(1-r\cos t)\sin t }{\pi (1-2r\cos t+r^2)}\right]^2(1-r\cos t)dt.$$

\n It is not clear if this integral can be simplified any further but it  is for sure related to the subject of elliptic integrals. For practical purposes it can be used to evaluate the probability numerically.
For instance, if $r=1/2$ one gets ${\cal P}\approx 0.1250$ and the graph of ${\cal P}_r$ as a function of $r\in [0,1]$ is included in Figure~12.

\[
\underset{Figure\ 12}{ \epsfig{file=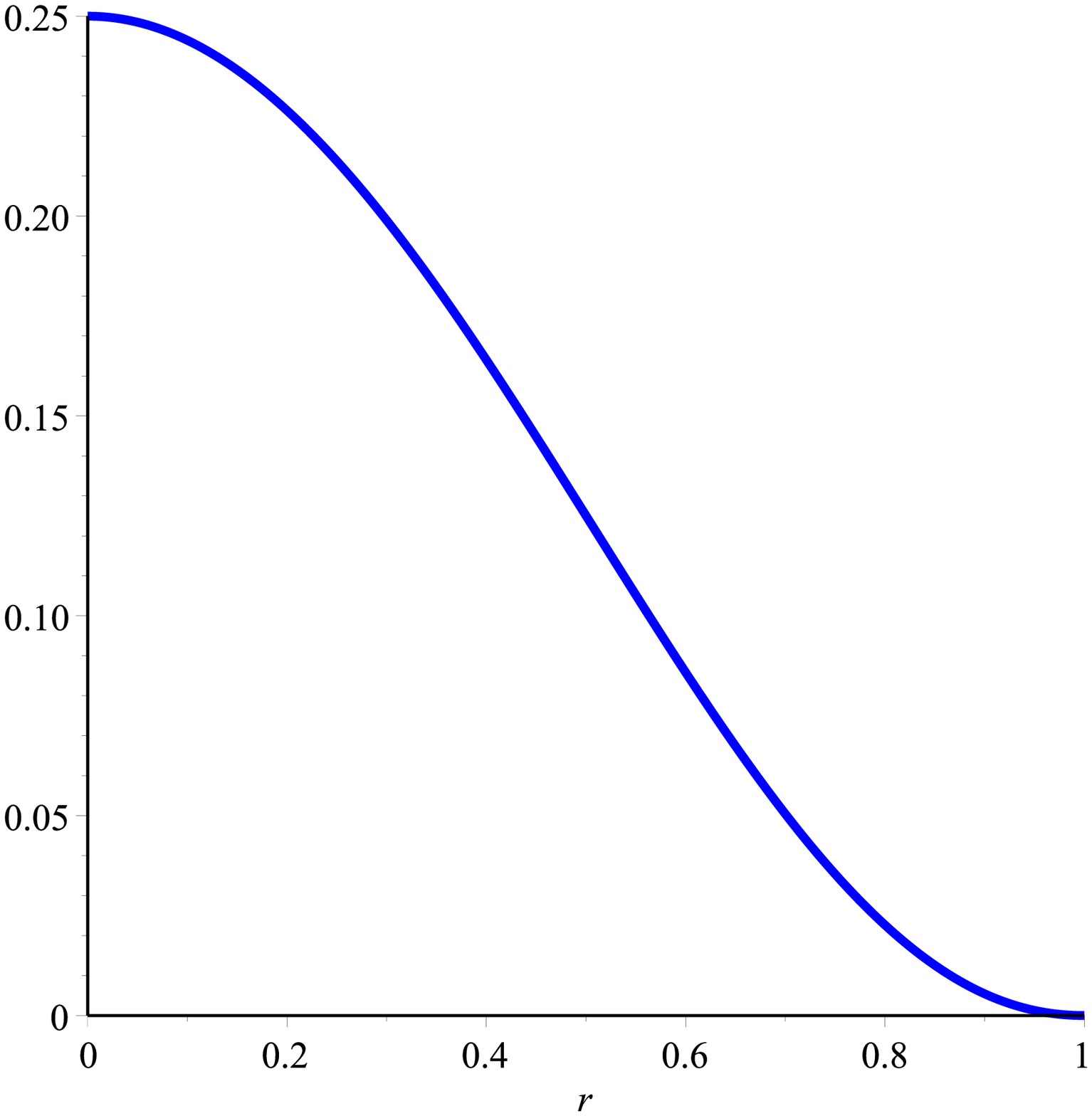,height=3in,width=3in}}
\]

\n The average value of ${\cal P}_r$ over  the disk is then

$$2\pi R^2 \int_0^1 {\cal P}_r rdr =0.07388002974=\frac{35}{48\pi^2}. $$
In other words, we recover the probability that four points chosen at random from the interior of a circle, the first point to be in the interior of the
triangle determined by the last three (see \cite{KM}).

\n

\end{document}